\documentclass[11pt,oneside]{article}
\usepackage{amssymb,amsmath}
\usepackage{url}
\usepackage[english]{babel}
\usepackage{pstricks,pst-node}
\renewcommand{\le}{\leqslant}
\renewcommand{\ge}{\geqslant}
\newcommand{\bad}{\mathbf{Bad}}
\renewcommand{\L}{{\rm L}}
\newcommand{\J}{{\rm J}}

\newcommand{\area}{\mathbf{area}}

\newcommand{\RR}{\mathbb{R}}
\newcommand{\ZZ}{\mathbb{Z}}
\newcommand{\QQ}{\mathbb{Q}}
\newcommand{\KK}{\mathbf{K}}
\newcommand{\MM}{\mathcal{M}}
\newcommand{\NN}{\mathbb{N}}

\newcommand{\JJ}{\mathcal{J}}
\newcommand{\II}{\mathcal{I}}
\newcommand{\LLL}{\mathcal{L}}
\newcommand{\KKK}{\mathbf{K}}
\newcommand{\RRR}{\mathcal{R}}
\newcommand{\TTT}{\mathcal{T}}
\newcommand{\vz}{\mathbf{z}}
\newcommand{\vw}{\mathbf{w}}

\newcommand{\vu}{\mathbf{u}}
\newcommand{\vv}{\mathbf{v}}

\newcommand{\vq}{\mathbf{q}}
\newcommand{\vp}{\mathbf{p}}

\newtheorem{lemma}{Lemma}
\newtheorem{theorem}{Theorem}
\newtheorem{proposition}{Proposition}

\newtheorem{SC}{Theorem \ref{sc}$^\prime$}

\newtheorem{corollary}{Corollary}

\newtheorem{vsec}{Proposition}

\newcommand{\cC}{{\cal C}}

\newcommand{\cJ}{{\cal J}}

\begin{document}

\title{On a problem in simultaneous Diophantine approximation: Schmidt's conjecture}

\author{Dzmitry Badziahin \footnote{Research supported by EPSRC grant EP/E061613/1}
\\ {\small\sc York}
\and  Andrew Pollington\footnote{Research supported by the National
Science Foundation}
\\ {\small\sc Arlington, VA}
  \and  Sanju Velani\footnote{Research supported by EPSRC grants EP/E061613/1 and EP/F027028/1 }
\\ {\small\sc York}}

\date{\it Dedicated to Graham Everest}

\maketitle

\begin{abstract}
For any $i,j \ge 0$ with $i+j =1$, let  $\bad(i,j)$ denote the set
of  points $(x,y) \in \RR^2$ for which  $ \max \{ \|qx\|^{1/i}, \,
\|qy\|^{1/j} \} > c/q $ for all $ q  \in \NN $. Here $c = c(x,y)$ is
a positive constant. Our main result implies that any finite
intersection of such sets has full dimension. This settles a
conjecture of Wolfgang M. Schmidt in the theory of simultaneous
Diophantine approximation.
%
%
%For any $i,j \ge 0$ with $i+j =1$, let  $\bad(i,j)$ denote the set
%of  points $(x,y) \in \RR^2$ for which there exists a positive constant $c(x,y)$  such that
%$$ \max \{ \|qx\|^{1/i}, \,
%\|qy\|^{1/j} \} > c(x,y)/q  \quad  \forall \ q  \in \NN    \, .$$
%Our main result implies that any finite
%intersection of such sets has full dimension. This settles a
%conjecture of Wolfgang Schmidt in the theory of simultaneous
%Diophantine approximation.

%This settles (a stronger form of) a conjecture of W.M. Schmidt in
%the theory of simultaneous Diophantine approximation.
\end{abstract}

%SV's draft

\section{Introduction}

A real number $x$ is said to be {\em badly approximable} if there
exists a positive constant $c(x)$ such that
\begin{equation*}
 \|qx\| \ > \ c(x) \ q^{-1}  \quad \forall \  q \in \NN   \ .
\end{equation*}
Here and throughout $ \| \cdot  \| $ denotes the distance of a real
number to the nearest integer. It is well know that the set $\bad$
of badly approximable numbers is of Lebesgue measure zero. However,
a result of Jarn\'{\i}k (1928) states that
\begin{equation}
\label{badfulldim} \dim \bad = 1  \ ,
\end{equation}
where $\dim X $ denotes  the Hausdorff dimension of the set $X$.
Thus, in terms of dimension the set of badly
approximable numbers is maximal; it has the same dimension as the
real line.  For details regarding  Hausdorff  dimension the reader is referred to  \cite{falc}.

In higher dimensions there are various natural generalizations of
$\bad$. Restricting our attention to the plane $\RR^2$, given a pair
of real numbers $i$ and $j$ such that
\begin{equation}\label{neq1}
0\le i,j\le 1   \quad {\rm and \  } \quad  i+j=1  \, ,
\end{equation}
a point $(x,y) \in \RR^2$ is said to be {\em $(i,j)$-badly
approximable} if there exists a positive constant $c(x,y)$ such that
\begin{equation*}
  \max \{ \; \|qx\|^{1/i} \; , \ \|qy\|^{1/j} \,  \} \ > \
  c(x,y) \ q^{-1} \quad \forall \  q \in \NN   \ .
\end{equation*}
Denote by   $\bad(i,j)$  the set of $(i,j)$-badly approximable
points in $\RR^2$.  If $i=0$, then we use the convention that $ x^{1/i}\::=0 $ and so $\bad(0,1)$ is identified with
$\RR \times \bad $. That is, $\bad(0,1)$ consists of points $(x,y)$
with $x \in \RR$ and $y \in \bad$.  The roles of $x$ and $y$ are
reversed if $j=0$. It easily follows from classical results in the
theory of metric Diophantine approximation   that $\bad(i,j)$ is of
(two-dimensional) Lebesgue measure zero. Building upon the work of
Davenport \cite{dav} from 1964, it has recently been shown in
\cite{PV} that $\dim \bad(i,j)=2$.  For further background and
various strengthenings of this full dimension statement the reader
is referred to \cite{KW,KTV,PV}. A consequence of the main result
obtained in this paper is the following   statement.

\begin{theorem}\label{sc}
Let  $(i_1,j_1), \ldots ,(i_d,j_d)$ be a finite number of pairs of
real numbers  satisfying~\eqref{neq1}. Then
%Let  $(i_t,j_t)$ be a countable number of pairs satisfying
%\eqref{neq1}. Then
$$
\dim \Big(\bigcap_{t=1}^{d} \bad(i_t,j_t) \Big) = 2   \ .
$$
\end{theorem}

\noindent Thus, the  intersection of any finitely  many badly
approximable sets $\bad(i,j)$ is  trivially non-empty and thereby
establishes the following conjecture of Wolfgang M. Schmidt \cite{schconj}
from the eighties.

\vspace{3ex}

\noindent{\bf Schmidt's conjecture \ } {\em For any $(i_1,j_1)$ and
$(i_2,j_2)$ satisfying \eqref{neq1}, we have that
 $$ \bad(i_1,j_1) \cap \bad(i_2,j_2) \  \neq \  \emptyset \ .
 $$}

\noindent  To be precise,  Schmidt stated the specific problem with
$i_1=1/3$ and $j_1=2/3$  and even this has previously resisted attack.  Indeed, the
statement
$$
\dim (\bad(1,0) \cap \bad(0,1) \cap \bad(i,j) ) = 2
$$
first obtained in \cite{PV}  sums up all previously known results.

\smallskip

 As noted by Schmidt, a counterexample to his conjecture would
imply the famous Littlewood conjecture: for any $(x,y) \in \RR^2$
$$
\liminf_{q \to \infty} q \, \|qx\| \, \|qy\| = 0    \; . $$ Indeed,
the same conclusion is valid if there exists any  finite (or indeed countable) collection
of  pairs $(i_t,j_t)$ satisfying \eqref{neq1} for which the intersection of the sets $ \bad(i_t,j_t) $ is empty. However, Theorem
\ref{sc} implies that no such finite collection exists  and
Littlewood's conjecture remains very much  alive and kicking. For background
and recent developments regarding Littlewood's conjecture see
\cite{pvl,vent}.

\subsection{The main theorem}

The key to establishing Theorem \ref{sc} is to investigate the
intersection  of the sets $ \bad(i_t,j_t) $ along fixed vertical
lines in the $(x,y)$-plane. With this in mind,  let $\L_x$ denote
the line parallel to the $y$-axis passing through the point $(x,0)$.
Next, for any real number $0\le i \le 1$, define the
 set
$$
\bad(i) \, := \, \{ x \in \RR : \exists \ c(x) > 0 \ {\rm so \ that}
\ \ \|q x\|
 > c(x) \, q^{-1/i}  \ \ \forall \ q \in \NN \} \; .
$$
Clearly,
\begin{equation}
\label{badi1}
 \bad = \bad(1) \, \subseteq \,   \bad(i)  \, ,
\end{equation}
which together with (\ref{badfulldim}) implies that
\begin{equation}
\label{badi} \dim \bad(i) = 1  \quad \forall \  i \in [0,1] \  .
\end{equation}
In fact, a straightforward argument involving the Borel-Cantelli
lemma from probability theory enables  us to conclude that for  $i
<1$ the complement of  $\bad(i)$ is of Lebesgue measure zero.  In
other words, for  $i <1$ the set $\bad(i)$ is not only of full
dimension but of full measure.

We are now in the position to state our main theorem.

\begin{theorem}\label{tsc}
Let  $(i_t,j_t)$ be a countable number of pairs of real numbers
satisfying \eqref{neq1} and let $i := \sup\{i_t : t \in \NN \}$.
Suppose that
\begin{equation} \label{liminfassump}
\liminf_{t \to \infty }  \min\{ i_t,j_t  \}  > 0   \ .
\end{equation}
Then, for any $\theta \in \bad(i)$ we have that
$$
\dim \Big(\bigcap_{t=1}^{\infty} \bad(i_t,j_t)   \cap \L_{\theta}
\Big) = 1 \ .
$$
\end{theorem}

\vspace*{2ex}

The hypothesis  imposed on $\theta $ is absolutely
necessary. Indeed, for $\theta  \notin \bad(i)$ it is readily
verified that the intersection of the sets $\bad(i_t,j_t)$ along the
line $ \L_{\theta}$  is empty --  see \S\ref{dfs} for the details.
However, in view of (\ref{badi1}), the dependence of $\theta$  on $i$ and
therefore the pairs $(i_t,j_t)$ can be entirely removed by insisting
that $ \theta \in \bad$. Obviously, the resulting statement is
cleaner but nevertheless  weaker than Theorem \ref{tsc}.

On the other hand, the  statement of Theorem \ref{tsc} is almost certainly valid without imposing the `$\liminf$' condition.   Indeed, this is trivially true  if the    number of $(i_t,j_t)$ pairs is finite.  In the course of establishing the theorem, it will become evident that in the countable `infinite' case   we require  \eqref{liminfassump} for an important but nevertheless technical reason.  It would be desirable to remove   \eqref{liminfassump}  from the statement of the theorem.

%On the other hand, the  statement of Theorem is almost certainly valid without imposing the `$\liminf$' condition.   Indeed, this is trivially true  if the    number of $(i_t,j_t)$ pairs is finite.  As far as we are concerned, in establishing the theorem in the countable `infinite' case   we require  \eqref{liminfassump} for an important but nevertheless technical reason.  It would be highly desirable to remove   \eqref{liminfassump}  from the statement of the theorem.

\vspace*{1ex}

The following corollary is technically far easier to establish than
the theorem and is more than adequate for establishing Schmidt's
conjecture.

\begin{corollary}\label{tscc}
Let  $(i_1,j_1), \ldots ,(i_d,j_d)$ be a finite number of pairs of
real numbers  satisfying~\eqref{neq1}. Then, for any $\theta \in
\bad$ we have that $$ \bigcap_{t=1}^{d} \bad(i_t,j_t)   \cap
\L_{\theta} \neq \emptyset \ . $$
\end{corollary}

\noindent We  give a self contained proof of the corollary during
the course of establishing Theorem \ref{tsc}.

\vspace*{2ex}

\noindent{\em Remark.  } The corollary is of independent interest
even when $d=1$. Since the work of  Davenport  \cite{dav}, it has
been known that there exist badly approximable numbers $x$ and $y$
such that $(x,y)$ is also a  badly approximable pair; i.e.
$\bad(1,0) \cap \bad(0,1) \cap \bad(1/2,1/2) \neq \emptyset$.
However it was not possible, using  previous methods, to specify
which $x$ one might take.  Corollary \ref{tscc} implies that we can
take $x$ to be any badly approximable number. So, for example, there
exist $y \in \bad$ such that $(\sqrt{2},y)\in\bad(1/2,1/2)$.
Moreover, Theorem \ref{tsc} implies that
$$
\dim \big( \, \{  y \in \bad : (\sqrt{2},y)\in\bad(1/2,1/2) \, \} \;
\big)   = 1  \ .
$$

\vspace*{.5ex}

\subsection{ Theorem \ref{tsc} $ \  \Longrightarrow \ $  Theorem \ref{sc} }
%It is appropriate to end this section by returning to Theorem \ref{sc}.

We show that Theorem \ref{tsc} implies the  following countable version of Theorem \ref{sc}.

\begin{SC} \label{scsv}
Let  $(i_t,j_t)$ be a countable number of pairs of real numbers
satisfying \eqref{neq1}. Suppose that
\eqref{liminfassump} is also satisfied.
%\liminf_{t \to \infty }  \min\{ i_t,j_t  \}  > 0   \ .
%\end{equation*}
Then
$$
\dim \Big(\bigcap_{t=1}^{\infty} \bad(i_t,j_t) \Big) = 2   \ .
$$
\end{SC}

\noindent Note that if the  number of $(i_t,j_t)$ pairs is finite, the `$\liminf$' condition is trivially satisfied and Theorem~\ref{sc}$^\prime$   reduces to Theorem \ref{sc}.

We proceed to establish Theorem~\ref{sc}$^\prime$   modulo Theorem \ref{tsc}.  Since any  set $\bad(i,j) $ is a subset of $\RR^2$, we immediately
obtain the upper bound result that
\begin{equation*}
 \dim \Big(\bigcap_{t=1}^{\infty} \bad(i_t,j_t)
\Big)  \le 2 \ .
\end{equation*}

\noindent The following general result that  relates the dimension
of a set to the dimensions of parallel sections, enables us to
establish the complementary lower bound estimate  -- see  \cite[pg.
99]{falc}.

%In view of (\ref{badi}) and  Theorem \ref{tsc},  the following
%general statement   enables us to establish the complementary
%lower bound estimate  -- see  \cite[pg. 99]{falc}.

\begin{vsec}
 Let $F$ be a subset of $\RR^2$ and let $E$ be
a subset of the $x$-axis. If $\dim (F \cap \L_x)  \ge t $ for all
$x \in E$, then $\dim  F  \ge  t + \dim E$.
\end{vsec}

\noindent With reference to the proposition, let $ F $ be a
countable intersection of $\bad(i,j) $ sets and let $E $ be the set
$ \bad $. In view of (\ref{badfulldim}) and  Theorem \ref{tsc}, the
lower bound result immediately follows.  Since (\ref{badfulldim}) is
classical and the upper bound statement for the dimension is
trivial, the main ingredient in establishing
Theorem~\ref{sc}$^\prime$ (and therefore Theorem \ref{sc})  is
Theorem \ref{tsc}.

\vspace*{2ex}

\noindent{\em Remark.  }
It is self evident that removing  \eqref{liminfassump} from the statement of Theorem \ref{tsc} would enable us to remove  \eqref{liminfassump} from the statement of Theorem~\ref{sc}$^\prime$. In other words, it would enable us to established  in full the  countable version of Schmidt's conjecture.

\subsection{The dual form  \label{dfs} }

At the heart of the proof of Theorem \ref{tsc}  is an intervals
construction that enables us to conclude that
$$
\bad(i,j)   \cap \L_{\theta}  \ \neq \ \emptyset  \qquad  \forall \
\theta \in \bad(i) \ .
$$
Note that this is essentially  the statement of Corollary  \ref{tscc} with  $d=1$.
The case when either $i=0$ or $j=0$ is relatively
straightforward so let us assume that
\begin{equation}\label{neq1sv}
0  <  i,j   <  1   \quad {\rm and \  } \quad  i+j=1  \, .
\end{equation}
In order to carry out the construction alluded to above, we shall work with the
equivalent dual form representation of the set $ \bad(i,j)$. % and for
%convenience we stay within  the unit square $[0,1)^2$.
In other words, a point $(x,y)\in  \bad(i,j)$ if there exists a
positive constant $c(x,y)$ such that
\begin{equation}\label{badij2}
\max\{|A|^{1/i}, |B|^{1/j}\}   \; \|Ax-By\|  >  c(x,y)  \qquad
\forall  \ (A,B)  \in\ZZ^2 \backslash \{(0,0)\}  \  .
\end{equation}

\noindent Consider for the moment the case $B=0$. Then,
(\ref{badij2}) simplifies to the statement that
$$ |A|^{1/i}  \; \|Ax\| >  c(x)  \qquad \forall  \ A  \in\ZZ\backslash\{0\}  \  .$$
It now becomes obvious that for a point $(x,y)$ in the plane  to
have any chance of  being   in $\bad(i,j)$,  we must have that $x
\in \bad(i)$. Otherwise, (\ref{badij2}) is violated and $\bad(i,j)
\cap \L_{x} = \emptyset $. This justifies the hypothesis imposed  on
$\theta$ in Theorem \ref{tsc}.

%
%Thus, for points $(\theta,y)$ on the vertical line $\L_{\theta}$
%to be in $\bad(i,j)$ we must have that $\theta  \in \bad(i)$.
%Otherwise, $\bad(i,j) \cap \L_{\theta} = \emptyset $.   This
%justifies the hypothesis imposed  on $\theta$ in Theorem
%\ref{tsc}.

\medskip

For $i$ and $j $ satisfying \eqref{neq1sv}, the equivalence of the `simultaneous' and `dual'  forms of $
\bad(i,j)$ is a consequence of the transference principle described
in \cite[Chapter 5]{cassels}. To be absolutely precise, without
obvious modification, the principle as stated in \cite{cassels} only
implies the equivalence in the case $i=j=1/2$. In view of this and
for the sake of completeness, we have included the modified
statement and its proof as an appendix.

\vspace{2ex}

\noindent{\em Notation.\ \ }For a real number $r$ we denote by $[r]$
its integer part and by $\lceil r\rceil$ the smallest integer
not less than $r$. For a subset $X$ of $\RR^n $ we denote by $|X|$
its Lebesgue measure.

\section{The overall strategy \label{strategy} }

Fix  $i$ and $j $ satisfying \eqref{neq1sv}   and  $\theta \in
\bad(i)$ satisfying  $ 0 < \theta < 1$.  Let $\Theta$ denote the
segment of the vertical line $\L_{\theta}$ lying within the unit
square; i.e.
$$\Theta:=  \{(x,y)\;:\;x=\theta,y\in [0,1]\}   \ . $$
In the section we describe the basic intervals construction that
enables us to conclude that
$$
\bad(i,j)   \cap \Theta  \ \neq \ \emptyset   \ .
$$

\vspace*{2ex}

%\noindent{\em Remark.  }
As  mentioned in  \S\ref{dfs}, the basic construction  lies at the heart of establishing  Theorem~\ref{tsc}.

%This clearly implies the statement of Corollary}  \ref{tscc} with $d=1$.

\vspace*{2ex}

\subsection{The sets $\bad_{c}(i,j)$  \label{secbadc}}  For any
constant $c >0$, let  $\bad_{c}(i,j)$  denote  the set of points
$(x,y)\in \RR^2$  such that
\begin{equation}\label{neq2}
\max\{|A|^{1/i}, |B|^{1/j}\}   \; \|Ax-By\|  >  c   \qquad \forall \
(A,B)  \in\ZZ^2 \backslash \{(0,0)\}  \  .
\end{equation}

%
%\max\{|A|^{1/i}, |B|^{1/j}\}|Ax-By+C|\ge c    \qquad
%  \ \forall (A,B,C)\in\ZZ^3 {\rm \ with \ } (A,B)\neq (0,0)\;.

\noindent It is easily seen that $  \bad_{c}(i,j) \subset \bad(i,j)$
and
 $$\bad (i,j)    \, =   \, \bigcup_{c > 0}  \bad_{c}(i,j)    \ . $$
Geometrically,  given  integers  $A,B,C$  with $(A,B)\neq (0,0)$
consider the line $L=L(A,B,C)$ defined by the equation
$$
Ax-By+C = 0   \ .
$$
The set $\bad_{c}(i,j)$ simply consists of points in the plane that
avoid  the  $$ \frac{c}{ \max\{|A|^{1/i}, |B|^{1/j}\} } $$
thickening  of each line $L$  -- alternatively, points in the plane
that lie within any such neighbourhood are removed.   With reference
to our fixed $\theta \in \bad(i)$, let us assume that
\begin{equation} \label{svs}
c(\theta) \ge c > 0  \ .
\end{equation}
Then, by definition
\begin{equation}\label{neq3}
|A|^{1/i}  \; \|A  \theta \| >   c   \qquad \forall  \ A
\in\ZZ\backslash\{0\}  \
\end{equation}
and the  line $L_\theta$ (and  therefore the segment $\Theta$)  avoids
the thickening of any vertical line $L=L(A,0,C)$.  \emph{Thus,
without loss of generality, we can assume that $B\neq 0$.} With this
in mind, it is easily verified that the thickening of a line
$L=L(A,B,C)$ will remove from $\Theta$ an interval  $\Delta(L)$
centered at $(\theta,y)$ with
$$
y=\frac{A\theta+C}{B}
$$
and length
%$$
%|\Delta(L)| = \frac{2c}{H(A,B)}
%$$
%where
%$$
%    H(A,B)  := |B|\max\{|A|^{1/i},|B|^{1/j}\}   \, .
%$$
%
\begin{equation}\label{height}
|\Delta(L)| = \frac{2c}{H(A,B)}  \qquad   {\rm where  \ } \qquad
    H(A,B)  := |B|\max\{|A|^{1/i},|B|^{1/j}\}   \, .
\end{equation}

\noindent For reasons that will soon become apparent, the quantity
$H(A,B)$   will be referred to as the \emph{height of the line}
$L(A,B,C)$.  In short, the height determines the amount of material
a line  removes from the fixed vertical line $\L_{\theta}$ and
therefore from  $\Theta$.

 The upshot of the above  analysis is that the set
$$
\bad_c(i,j)   \cap \Theta
$$
consists of points $(\theta, y)$ in the unit square which avoid all
intervals  $\Delta(L)$ arising from lines $L=L(A,B,C)$ with $B \neq
0$. Since
$$
\bad_c(i,j)   \cap \Theta \subset \bad (i,j)   \cap \Theta  \ ,
$$
the name of the game  is to show that we have something left after
removing these intervals.

\vspace*{1ex}

\noindent{\em Remark 1. \, } The fact that we have restricted our
attention to $\Theta$ rather than working on the whole line
$\L_{\theta}$ is mainly for convenience.  It also means that for any
fixed $A$ and $B$,  there are only a  finite number of lines
$L=L(A,B,C)$ of interest; i.e. lines for which $ \Delta(L) \cap
\Theta \neq \emptyset $. Indeed, with $c \le 1/2$ the number of
such lines is bounded above by $|B| +  2$.

\vspace*{1ex}

\noindent{\em Remark 2.  } Without loss of generality,  when
considering lines $L=L(A,B,C)$ we will~assume that
\begin{equation}\label{heighttriv}
(A,B,C)=1  \quad {\rm  and}  \quad  B>0   \; .
\end{equation}
 Otherwise we can divide  the
coefficients of $L$  by their common divisor  or by $-1$. Then the
resulting line $L'$ will satisfy the required conditions and
moreover $\Delta(L') \supseteq \Delta(L)$.  Therefore, removing the
interval $\Delta(L')$ from $\Theta$ takes care of removing
$\Delta(L)$.

\vspace*{1ex}

\noindent Note that  in view of (\ref{heighttriv}), for any line
$L=L(A,B,C)$ we always have that $ H(A,B)  \ge 1 $.

\vspace*{2ex}

\subsection{Description of basic construction}\label{secdescrip}

Let $R \ge  2$ be an integer.
%For convenience, we assume that $R$ is a power of two. Thus, $R=2^r $ for some integer  $r$.
Choose $c_1=c_1(R)$ sufficiently small
so that
\begin{equation} \label{eqc1}
c_1   \ \le \ \textstyle{\frac{1}{4}}   R^{-\frac{3i}{j}}
\end{equation}
and
\begin{equation} \label{eq_c1}
c \ := \   \frac{c_1}{R^{1+\alpha}}
\end{equation}
satisfies \eqref{svs} with
\begin{equation} \label{eqc1sv}
\alpha   \, :=  \, \textstyle{\frac{1}{4}}  \,  ij    \,  .
\end{equation}

%Let $R \ge  2$ be an integer.
%%For convenience, we assume that $R$ is a power of two. Thus, $R=2^r $ for some integer  $r$.
%Choose $c=c(R)$ sufficiently small
%so that \eqref{svs} is satisfied and
%\begin{equation} \label{eqc1}
%c_1 := cR^{1+\alpha}  \ \le \ \textstyle{\frac{1}{4}}   R^{-\frac{3i}{j}}
%\end{equation}
%where
%\begin{equation} \label{eqc1sv}
%\alpha   \, :=  \, \frac{ij}{4}   \qquad (i,j \neq 0) \,  .
%\end{equation}

We  now  describe the basic construction that enables us to conclude
that
\begin{equation} \label{dag}
\bad_c(i,j)   \cap \Theta     \neq  \emptyset      \ .
\end{equation}
We start by subdividing  the segment  $\Theta$ from the $(\theta,0)$ end into
closed intervals  $J_0$ of equal length $c_1$.  Denote by $\JJ_0$
the collection of intervals $J_0$.
Thus,
$$
\# \JJ_0    =  [c_1^{-1}]  \ .
$$
 The idea is to establish,  by induction on  $n$, the
existence of a collection $\JJ_n$ of closed intervals $J_n$ such
that $\JJ_n $ is nested in $\JJ_{n-1}$;   that is, each interval
$J_n$ in $\JJ_n$ is contained in some interval $J_{n-1}$ in
$\JJ_{n-1}$.    The length of an interval  $J_n$  will be given by
$$ |J_n|  \, :=  \,   c_1    \, R^{-n}   \  ,
$$
and each  interval $J_n $ in $\JJ_n$ will satisfy the condition that
\begin{equation}  \label{cond}
J_{n} \,  \cap \, \Delta(L) \, =  \, \emptyset   \qquad    \forall \
\  L=L(A,B,C)  \ \  \mbox{with    } \   H(A,B)  <   R^{n-1}   \,   .
\end{equation}
In particular, we put
$$
\KK_{c}    =   \KK_{c(R)}:=       \bigcap_{n=1}^\infty
\bigcup_{J\in\JJ_n}J      \ .
$$
By construction, we have that
$$
\KK_c  \subset \bad_c(i,j)   \cap \Theta   \ .
$$
Moreover, since the intervals $J_n$  are nested, in order to
establish  \eqref{dag} it suffices to show that each  $\JJ_n$
is non-empty; i.e.
$$\# \JJ_n   \ge 1    \qquad    \forall  \  n = 0,1, \ldots    \  . $$

\vspace{4ex}

%\subsubsection{The induction}
\noindent{\bf  The induction. } For $n=0$, we trivially have that
(\ref{cond}) is satisfied for any interval $J_0\in \JJ_0$. The
point is that in view of (\ref{heighttriv}) there are no lines
satisfying the height condition  $H(A,B)<1$.
For the same reason (\ref{cond})  with $n=1$ is trivially satisfied for any interval $J_1$  obtained by subdividing each $J_0$ in
$\JJ_0$  into $R$ closed intervals of equal length $c_1 R^{-1} $.  Denote by $\JJ_1 $  the resulting collection of intervals $J_1$ and note that
$$
\# \JJ_1    =  [c_1^{-1}]  \, R  \ .
$$
In general, given $\JJ_n$
satisfying  (\ref{cond}) we wish to construct a nested
collection $\JJ_{n+1}$ of intervals $J_{n+1}$ for which
(\ref{cond}) is satisfied with $n$ replaced by $n+1$. By
definition, any interval $J_n$ in $\JJ_n$ avoids intervals
$\Delta(L)$ arising from lines with height bounded above by
$R^{n-1}$.  Since any `new' interval $J_{n+1}$ is to be nested in
some $J_n$,  it is enough to show that $J_{n+1}$ avoids intervals
$\Delta(L)$ arising from lines $L=L(A,B,C)$ with height satisfying
\begin{equation}\label{zeq2}
R^{n-1}\le H(A,B)<R^n  \   .
\end{equation}
Denote by $\cC(n)$ the collection of all lines satisfying this height
condition.  Throughout, we are already  assuming that lines satisfy
(\ref{heighttriv}). Thus,  formally
$$
\cC(n) := \left\{L= L(A,B,C)  \, : \,  L  \ \ {\rm satsifies \
(\ref{heighttriv})   {\rm \ and \ } (\ref{zeq2}) \, }   \right\}   \
$$
and it is precisely this  collection of lines that comes into play
when  constructing   $\JJ_{n+1}$ from $\JJ_{n}$.  We now proceed
with the construction.

%The upshot of the above is that it is this collection $\cC(n)$ of lines that needs to be considered when  constructing   $\JJ_{n+1}$ from $\JJ_{n}$.
%
%The upshot is that in order to  construct  $\JJ_{n+1}$  we need to naturally  brings into play lines
%from the collection $\cC(n)$.
%
%Suppose $J_n$ is a good interval; that is, all points $(\theta,y)$  avoid intervals  $\Delta(L)$ arising from lines with
%height less than $R^{n-1}$. In short, $J_n\in \JJ_n$.
%Suppose $J_n$ is a good interval; that is, $J_n\in \JJ_n$. We now
%construct a collection of good intervals $J_{n+1}$ lying within
%$J_n$ in two stages.

\vspace{2ex}

\noindent{\em  Stage 1: The collection $\II_{n+1}$. }  We subdivide each $J_n$ in
$\JJ_n$  into $R$ closed intervals $I_{n+1}$ of equal length and
denote by $\II_{n+1}$ the collection of such intervals. Thus,
$$|I_{n+1}|=c_1 R^{-n-1}       \qquad {\rm and \ }  \qquad   \#
\II_{n+1} =     R  \,  \times   \,  \# \JJ_n    \  .$$ In view of the nested
requirement, the collection $\JJ_{n+1}$ which  we are attempting to
construct will  be a  sub-collection  of  $\II_{n+1}$. In other
words, the intervals $I_{n+1}$ represent possible candidates for
$J_{n+1}$. The goal now is simple -- it is  to remove those `bad'
intervals $I_{n+1}$ from $\II_{n+1}$  for which
\begin{equation}  \label{svt}
I_{n+1} \,  \cap \, \Delta(L) \, \neq  \, \emptyset   \ \ \mbox{ for
some \ } L \in \cC(n)  \ .
\end{equation}
%Trivially, the collection $C(O)$ is empty and so  we simply put
%$$
%\JJ_1 \, := \, \II_1 \ .
%$$
Note that the number of bad intervals
that can be removed by any single line $L=L(A,B,C) $   is bounded by
\begin{equation}\label{neq4}
\frac{|\Delta(L)|}{|I_{n+1}|}  \, + \, 2  \; =  \;   2\frac{cR^{n+1}}{c_1H(A,B)} \, + \, 2 \; = \; \frac{2R^{n-\alpha}}{H(A,B)}  \, +  \, 2   \ .
\end{equation}
Thus any single line  $L$  in  $\cC(n) $ can remove up to $[2R^{1-\alpha}]+2$  intervals   from   $\II_{n+1}$.  Suppose,  we crudely  remove this maximum number  for each $L$ in  $\cC(n) $. Then, for $n$ large enough, a straightforward calculation shows  that  all the intervals from  $\II_{n+1}$  are eventually  removed  and the construction comes to a halt.  In other words, we need to be much more sophisticated in our approach.

\vspace*{2ex}

\noindent{\em  Stage 2: Trimming. }~Even before considering
the effect that lines from $\cC(n)$ have on intervals in $\II_{n+1}$,
 we trim the collection $\II_{n+1}$ by removing from each $J_n$
the first $\lceil R^{1-\alpha}\rceil$ sub-intervals $I_{n+1}$ from each end.
Let us denote by $J_n^-$ the resulting `trimmed' interval and by
$\II_{n+1}^{\, -} $ the resulting `trimmed' collection. This process
removes $\# \JJ_n \times 2 \, \lceil R^{1-\alpha}\rceil$ intervals
$I_{n+1}$  from $\II_{n+1}$ regardless of  whether an interval is
bad or not. However, it  ensures that for any remaining interval
$I_{n+1} $ in  $\II_{n+1}^{\, -} $  which satisfies  (\ref{svt}) the
line $L$ itself must intersect the associated interval $J_n$ within
which $I_{n+1} $  is nested.   The upshot of `trimming' is that when
considering (\ref{svt}),  we  only need to consider those lines $L$
from $\cC(n)$ for which
$$
J_{n} \,  \cap \,L \, \neq  \, \emptyset   \ \ \mbox{ for some \ }
J_n  \in \JJ_n  \ .
$$
The intervals $\Delta(L)$  arising from the `other' lines  are either removed by the
trimming process or  they do
not even intersect intervals in $\JJ_n$ and therefore they  can not  possibly
remove any intervals from $\II_{n+1}$.

\vspace*{1ex}

The sought after collection $\JJ_{n+1}$ is precisely that
obtained by  removing those `bad' intervals $I_{n+1}$  from $\II_{n+1}^{\, -} $   which  satisfy (\ref{svt}). Formally, for $n \ge 1 $  we let
\begin{equation}\label{svgood}
\JJ_{n+1}   \,  : =  \,   \{I_{n+1}\in  \II_{n+1}^{\, -}  \;:\;    \Delta(L)\cap I_{n+1}  = \emptyset        \  \ \forall  \ \ L \in \cC(n)  \}    \ .
\end{equation}
For any   strictly positive $\epsilon <   \frac12  \, \alpha^2 $  and $R > R_0(\epsilon) $ sufficiently large, we claim that
\begin{equation} \label{countgood}
\# J_{n+1} \ \ge \  (R-5R^{1-\epsilon})  \, \times   \,   \# J_{n}    \qquad    \forall  \  n = 0,1, \ldots    \  .
\end{equation}
Clearly, this implies that
\begin{equation*} \label{countgoodtriv}
\# J_{n+1} \ \ge \  (R-5R^{1-\epsilon})^{n+1} \, >  \, 1 \,
\end{equation*}
which  in turn completes the proof of the induction step and therefore  establishes  $\eqref{dag}$. Thus, our goal now is to  justify (\ref{countgood}).

\vspace{4ex}

\noindent{\em  Stage 3: The sub-collection $\cC(n,l)$. }    In the
first instance we subdivide the collection $\cC(n)$ of lines  into various  sub-collections
that  reflect a common  geometric  configuration.  For any integer $ l \ge 0$, let $\cC(n,l) \subset \cC(n)  $ denote the collection of lines $L=L(A,B,C)$  satisfying  the additional condition that
\begin{equation}\label{zeq1}
R^{-\lambda(l+1)}R^{\frac{nj}{j+1}}  \, \le \, B  \, < \,
R^{-\lambda l}R^{\frac{nj}{j+1}}       \ \,
\end{equation}
where $$ \lambda := 3/j   > 1   \ .
$$
%In view of (\ref{zeq2}), it follows  that
%$$ 0 \le l  \ \le  \ \frac{nj}{\lambda(j+1)}  \, . $$
%Also, since $\lambda>1$, we have that $l<n$.
%
Thus the $B$ variable associated with any line  in $\cC(n,l)$ is within a tight range governed   by~(\ref{zeq1}).
In view of \eqref{zeq2}, it follows that
$
B^{1+1/j}< R^n
$
and so
$
 1 %\stackrel{\eqref{heighttriv}}{\le}
 \le B<R^{\frac{nj}{j+1}}    \ .
 $
Therefore,
$$ 0 \ \le \  l  \ <  \ \frac{nj}{\lambda(j+1)}  \ < \ n  \, .  $$
%and since $\lambda>1$ we have that $l<n$.
A useful `algebraic' consequence of imposing (\ref{zeq1}) is that
\begin{equation}\label{zeq3slv}
H(A,B) \;  = \;  |B| \, |A|^{1/i}   \qquad  \forall \  \
L(A,B,C)  \,  \in  \,   \cC(n,l >0) %{\rm \ with  \  }  l > 0
\, .
\end{equation}
To see this,  suppose that the
$B^{1/j}$  term is the  maximum  term associated with $H(A,B)$.  Then, by
\eqref{zeq2} we have that
$$
\textstyle{ B\cdot B^{1/j}\ge R^{n-1} \  \Longrightarrow   \  B  \, \ge  \, R^{\frac{(n-1)j}{j+1}}  \ . }
$$
Thus,  by definition  $L(A,B,C) \in \cC(n,0)$. Moreover,  in view of
(\ref{zeq3slv}) and the definition of $\cC(n,l)$,  it follows that
%$$
% R^{\lambda l-1}\cdot
%R^{\frac{n}{j+1}}  \; <  \ |A|^{1/i}  \ <  \; R^{\lambda(l+1)}R^{\frac{n}{j+1}}
%$$
%and so
\begin{equation}\label{zeq3}
R^{(\lambda l-1)i}\cdot
R^{\frac{ni}{j+1}}  \; <  \ |A| \ <  \; R^{\lambda(l+1)i}R^{\frac{ni}{j+1}}  \qquad  \forall \  \
L(A,B,C)  \,  \in  \,   \cC(n,l >0)  \, .
\end{equation}
The upshot is that for $l > 0$,   both the  $A$ and $B$  variables  associated with lines  in $ \cC(n,l)$   are tightly controlled.   The above consequences of imposing  (\ref{zeq1}) are important but are out weighed by the significance of the following `geometric'  consequence.

\begin{theorem}\label{thmlem3}
All lines from $\cC(n,l)$ that intersect a fixed interval $J_{n-l} \in \JJ_{n-l}$
pass through a single rational point $P$.
\end{theorem}

\noindent The theorem is proved in \S\ref{triangsec}.  It implies that if we have three or more lines from $\cC(n,l)$  passing through any fixed  interval $J_{n-l}$, then the lines  can not possibly  enclose a  triangular region.  In short, triangles are not allowed.  The theorem represents  a crucial ingredient towards establishing
the following  counting statement. \emph{Let $ l \ge 0 $ and  $J_{n-l} \in \JJ_{n-l}$.
Then, for any  strictly positive $\epsilon <  \frac12 \alpha^2 $  and $R > R_0(\epsilon) $ sufficiently large, we have  that
%
%$$
%\#\{I_{n+1}\in \II_{n+1}\;:\;
% \exists \ L\in \cC(n,l) {\rm \ such \ that \ } J_{n-l}  \cap \Delta(L)\cap
%I_{n+1}\neq\emptyset\}      \ \le \ R^{1-\epsilon}
%$$
%where $\epsilon <  \alpha^2 $ is some positive constant. By taking
%together all intervals $J_{n-l}\in \JJ_{n-l}$ we immediately have
%that
%$$
%\#\{I_{n+1}\in \II_{n+1}\;:\; \exists \ L \in \cC(n,l)  {\rm \ such \
%that \ }   \Delta(L)\cap I_{n+1}\neq\emptyset\}  \ \le \
%R^{1-\epsilon}  \times \, \#\JJ_{n-l}.
%$$
\begin{equation} \label{countbad}
\#\{I_{n+1}\in  \II_{n+1}^{\, -}\;:\;
  J_{n-l}  \cap \Delta(L)\cap
I_{n+1}\neq\emptyset    {\rm \  \ for  \ some \   }   \ L\in \cC(n,l)  \}   \ \le \ R^{1-\epsilon}   \ .
\end{equation} }

Armed with this estimate it is reasonably straightforward to
establish (\ref{countgood}).  We use induction. For $n=0$, we have
that
$$
\#  \JJ_1  \, = \, R  \, \times  \, \#  \JJ_0 \
$$
and so (\ref{countgood}) is obviously true.  For $n \ge 1$, we suppose that
$$
\# J_{k+1} \ \ge \  (R-5R^{1-\epsilon})  \, \times   \,   \# J_{k}      \qquad    \forall  \  k = 0,1, \ldots, n-1    \
$$
and proceed to establish the statement for $k=n$.  In view of (\ref{countbad}),  we  have that the total number    of                                                                                                                                             intervals $I_{n+1}$  removed from  $\II_{n+1}^{\, -} $  by lines from $\cC(n,l)$ is bounded above by
$$
 R^{1-\epsilon}  \, \times  \, \#  \JJ_{n-l}    \ .
$$
It now follows that
\begin{eqnarray}\label{eq18sv}
~ \hspace*{-10ex}\#\{I_{n+1}\in \II_{n+1}^{\, -}  \;  :  \;  \Delta(L) \!\!\!\!  & \cap  & \!\!\!\!
I_{n+1}\neq\emptyset  {\rm \  \ for  \ some \   }   \ L\in \cC(n)  \} \nonumber  \\[2ex]   \ & \le & \
\sum_{l=0}^{n} R^{1-\epsilon } \,  \#  \JJ_{n-l}    \ \le  \  R^{1-\epsilon } \,  \#  \JJ_{n}  \ + \
R^{1-\epsilon } \sum_{l=1}^n    \#  \JJ_{n-l}    \ .
\end{eqnarray}
In view of the  induction hypothesis,  for $R$ sufficiently large we have that
$$
\sum_{l=1}^n    \#  \JJ_{n-l}  \ \le  \  \#  \JJ_{n}   \sum_{l=1}^\infty
(R-5R^{1-\epsilon})^{-l}  \ \le \  2 \;   \#  \JJ_{n}  \
$$
and so
\begin{equation}\label{eq18svbb}
{\rm  l.h.s. \ of \ }  (\ref{eq18sv})    \ \le \ 3 \; R^{1-\epsilon
} \,  \#  \JJ_{n}.
\end{equation}
Therefore,  for $R$ sufficiently large
\begin{eqnarray}\label{eq18sv3}
  \#  \JJ_{n+1 }  \  & = &   \  \# \II_{n+1}^{\, -}  \ - \  {\rm  l.h.s. \ of \ }  (\ref{eq18sv})   \nonumber \\[1ex]   \nonumber
 & \ge &   (R  -  2 \, \lceil R^{1-\alpha}\rceil ) \, \#   \JJ_{n}   \  - \ 3 \; R^{1-\epsilon } \,  \#  \JJ_{n}  \nonumber \\[1ex] \nonumber
& =  &  (R-5 R^{1-\epsilon})  \   \#  \JJ_{n}  \  .
\end{eqnarray}

\noindent This completes the induction step and therefore establishes (\ref{countgood}).
Thus, our goal now is to  justify (\ref{countbad}).

\vspace{4ex}

\noindent{\em  Stage 4: The sub-collection $\cC(n,l,k)$. } Clearly,
when attempting to establish (\ref{countbad}) we are only interested
in lines $L=L(A,B,C)$ in $\cC(n,l)$ which  remove intervals. In
other words,   $ \Delta(L)\cap I_{n+1}\neq\emptyset $ for   some
$I_{n+1} \in \II_{n+1}^{\, -} $.  Now the total number of intervals
that a line  $L$ can remove depends on the actual value of its
height. In the situation under consideration, the height satisfies
(\ref{zeq2}).   Therefore,  in view of (\ref{neq4}) the total number
of intervals $I_{n+1}$ removed by $L$ can vary anywhere  between
$1$ and  $[2R^{1-\alpha}]+2$.  In a nutshell, this variation is too
large to handle and we need to introduce a tighter control on the
height. For any integer $ k \ge 0$, let $\cC(n,l,k) \subset \cC(n,l)
$ denote the collection of lines $L=L(A,B,C)$  satisfying  the
additional condition that
\begin{equation}\label{eq1}
2^kR^{n-1}  \ \le \ H(A,B) \ < \ 2^{k+1}R^{n-1} \; .
\end{equation}
In view of \eqref{zeq2}, it follows that
\begin{equation}\label{eq1vbvb}
0  \ \le \  k  \ <  \ \frac{\log R}{\log 2}  \ .
\end{equation}
%The following counting result is at the heart of establishing  (\ref{countbad}) and indeed represents the key to unlocking Schmidt's conjecture.
The following counting result implies  (\ref{countbad}) and indeed represents the technical key to unlocking Schmidt's conjecture.

\begin{theorem}\label{maincount}
Let $ l, k  \ge 0 $ and  $J_{n-l} \in \JJ_{n-l}$.
Then, for any  strictly positive $\epsilon <  \alpha^2 $  and $R > R_0(\epsilon) $ sufficiently large, we have  that
\begin{equation} \label{countbadmain}
\#\{I_{n+1}\in  \II_{n+1}^{\, -}\;:\;
  J_{n-l}  \cap \Delta(L)\cap
I_{n+1}\neq\emptyset    {\rm \  \ for  \ some \   }   \ L\in \cC(n,l,k)  \}   \ \le \ R^{1-\epsilon}   \ .
\end{equation}
\end{theorem}

%\noindent Clearly, Theorem \ref{thmlem3} is applicable to the collection $\cC(n,l,k)$  and it used in the proof of this theorem that we make use of Theorem

\noindent Theorem \ref{maincount} is proved in \S\ref{sec6}.  It is in this proof  that we make use of Theorem \ref{thmlem3}. Note that the latter is applicable since $\cC(n,l,k) \subset \cC(n,l)  $. Also note that in view of the `trimming' process, when considering  \eqref{countbadmain}   we can assume that  $ J_{n-l}  \cap L \neq\emptyset $.    With Theorem  \ref{maincount} at our disposal, it follows that
for $R$ sufficiently large
$$
{\rm  l.h.s. \ of \ } (\ref{countbad})  \ \le \ \frac{\log R}{\log 2}  \; \times \; R^{1-\epsilon}  \ \le \  R^{1-\frac{1}{2} \epsilon } \ .
$$
This establishes  (\ref{countbad}) and  completes the description of the   basic construction.

\vspace*{4ex}

\noindent {\em Remark. \, }  We emphasize that from the onset of
this section we have  fixed  $i$ and $j $ satisfying \eqref{neq1sv}.
Thus  this condition on $i$ and $j $  is implicit within the
statements of Theorems \ref{thmlem3} and \ref{maincount}.

%To summarize,   modulo  Theorems \ref{thmlem3} and \ref{maincount}, the basic construction yields the statement of Corollary  \ref{tscc} with $d=1$.  We now show that with very little extra effort we are actually  in the position to prove the corollary in full and thereby establish Schmidt's conjecture.

%The dimension result however needs more work. We will remove some
%segments from $\JJ_n$ to construct another sequence of sets $\MM_n$
%in such a way that any segment $J_n\in \MM_n$ contains at least $R -
%2R^{1-\epsilon}$ subsegments $J_{n+1}\in \MM_n$ where $\epsilon$ is
%some positive number.
%
%The idea is to show that $\MM_n$ is nonempty. This result will be
%proved in subsection ????.
%
%Now we can apply the same procedure for finitely many pairs
%$(i_t,j_t)$ simultaneously where $1\le t\le d$. Denote by
%$\MM_n(i_t,j_t)$ the set $\MM_n$ which we get for a particular pair
%$(i_t,j_t)$. Define
%$$
%\MM_n^*:= \bigcap_{t=1}^d \MM_n(i_t,j_t).
%$$
%
%
%Since for each pair $(i_t,j_t)$ and each segment $J\in
%\MM_n(i_t,j_t)$ there is no more than $2R^{1-\epsilon}$ subsegments
%is removed then for each $J\in \MM_n^*$ there is no more than
%$2dR^{1-\epsilon}$ subsegments is removed. Therefore the amount of
%surviving subsegments is at least
%$$
%R-2dR^{1-\epsilon}.
%$$
%

\section{Proof of Corollary~\ref{tscc}: Modulo Theorems \ref{thmlem3} and \ref{maincount}}\label{sec7}

Modulo  Theorems \ref{thmlem3} and \ref{maincount}, the basic
construction of \S\ref{secdescrip}  yields the statement of
Corollary~\ref{tscc} for any single  $(i,j)$ pair satisfying
\eqref{neq1sv}.  We now show that with very little extra effort, we
can modify the  basic construction to simultaneously incorporate any
finite number of $(i,j)$ pairs satisfying \eqref{neq1}. In turn,
this  will prove  Corollary~\ref{tscc} in full and thereby establish
Schmidt's~conjecture.

%
%We now show that with very little extra effort, we can actually prove the  corollary in full and thereby establish Schmidt's conjecture.

%Consider the pairs $(i_t,j_t)$ satisfying \eqref{neq1}. Let's do the
%procedure for all pairs simultaneously. Note that the imposed
%constants $R_t, c_{1t}$ and $c_{3t}$ will depend on the pair
%$(i_t,j_t)$. However we can take
%$$
%R:=\max_1\le t\le d R_t;\quad c_1:=\min_1\le t\le d c_{1t}.
%$$
%Since al the conditions on $R$ are lower bounds, all the conditions
%on $c_1$ are upper bounds then the constructed $c_1$ and $R$ will
%work for all pairs $(i_t,j_t)$ simultaneously.

\subsection{Modifying the basic construction for finite pairs}\label{mbc}
To start with we suppose that the $d$ given pairs
$(i_1,j_1),\ldots,(i_d,j_d)$  in Corollary~\ref{tscc}
satisfy~\eqref{neq1sv}. Note that for each $t=1,\ldots,d$, the
height $H(A,B)$  of a given line $L=L(A,B,C) $ is dependent on the
pair $(i_t,j_t)$. In view of this and with reference to
\S\ref{strategy}, let us write $H_t(A,B)$  for $H(A,B)$,
$\Delta_t(L)$ for $ \Delta(L)$   and $\cC_t(n)$ for $\cC(n)$. With
this in mind, let $R \ge  2$ be an integer. Choose $c_1=c_1(R)$
sufficiently small so that
%\begin{equation} \label{eqc1fc}
%c_1   \ \le \ \textstyle{\frac{1}{4}}   R^{-\frac{3i_t}{j_t}}   \qquad \forall \  \ 1 \le t \le d \,
%\end{equation}
\begin{equation*}
c_1   \ \le \ \textstyle{\frac{1}{4}}   R^{-3i_t/j_t}  \quad  \qquad \forall \  \  \ 1 \le t \le d \, ,
\end{equation*}
and for each $t=1,\ldots,d$
$$
c(t) \ := \   \frac{c_1}{R^{1+\alpha_t}}
$$
satisfies \eqref{svs} with
\begin{equation*}
\alpha_t  \, :=  \,   \textstyle{\frac{1}{4}}  \, i_t \, j_t    \,  .
\end{equation*}

\noindent Note that with this choice of $c_1$ we are able to separately carry out the basic construction of \S\ref{secdescrip} for each $(i_t,j_t)$ pair and therefore conclude that
$$
 \bad_{c(t)} (i_t,j_t)   \cap \Theta     \neq  \emptyset     \quad  \qquad \forall \  \  \ 1 \le t \le d \, .
$$
We now describe the minor modifications to  the  basic construction
that  enable us to simultaneously  deal with the $d$  given
$(i_t,j_t)$  pairs and  therefore conclude that
$$
\bigcap_{t=1}^d \bad_{c(t)} (i_t,j_t)   \cap \Theta     \neq  \emptyset      \ .
$$
The modifications are essentially at the `trimming' stage and in the manner in which  the collections $\JJ_{n}$ for $n \ge 2$ are defined.

Let $c_1$ be as above.  Define the collections  $\JJ_0$ and $\JJ_1$
as in the basic construction.   Also  Stage~1 of the `induction' in
which the collection $\II_{n+1}$ is introduced remains unchanged.
However, the goal now is to remove those `bad' intervals $I_{n+1}$
from $\II_{n+1}$ for which
\begin{equation}\label{bad_i}
I_{n+1}\cap \Delta_t(L)\neq \emptyset  \ \  \mbox{ for some } t=1,\ldots,d\mbox{ and }L\in \cC_t(n)   \; .
\end{equation}
Regarding Stage~2,  we trim  the collection $\II_{n+1}$ by removing from each $J_n$ the first $\lceil R^{1-\alpha_{\min}}\rceil$ sub-intervals $I_{n+1}$ from each end. Here
$$\alpha_{\min}  \, :=   \, \min\{\alpha_1,\ldots,\alpha_d\}  \; .$$
This gives rise to the trimmed collection $\II_{n+1}^-$ and we define $\JJ_{n+1}$ to be the collection obtained by removing those `bad' intervals $I_{n+1}$ from $\II_{n+1}^-$ which satisfy \eqref{bad_i}.  In other words, for $n \ge  1$ we let
$$
\JJ_{n+1}   \,  : =  \,   \{I_{n+1}\in  \II_{n+1}^{\, -}  \;:\;    \Delta_t(L)\cap I_{n+1}  = \emptyset        \quad  \forall  \ 1\le t\le d\  \mbox{ and }\ L \in \cC_t(n)  \}    \ .
$$

\noindent Apart from obvious notational modifications, Stages 3 and 4 remain pretty much unchanged and enable us to establish \eqref{eq18svbb} for each  $t=1,\ldots,d$. That is, for any  strictly positive $\epsilon <  \frac12 \alpha_t^2 $  and $R > R_0(\epsilon) $ sufficiently large
\begin{eqnarray}\label{eq18svbbdd}
~ \hspace*{-10ex}\#\{I_{n+1}\in \II_{n+1}^{\, -}  \;  :  \;  \Delta_t(L) \ \cap  \
I_{n+1}\neq\emptyset \   {\rm \  \ for  \ some \   } \!\!\!\!\!\!\!\!  &  & \!\!\!\!\!\!\!\!   \ L\in \cC_t(n)  \} \nonumber  \\[2ex]   \ & \le &
 3R^{1-\epsilon } \,  \#  \JJ_{n}    \ .
\end{eqnarray}
It follows that for any strictly positive $\epsilon <  \frac12  \, \alpha_{\min}^2 $
and $R > R_0(\epsilon) $ sufficiently large
%$$
%\# J_{n+1} \ \ge \  (R-5dR^{1-\epsilon})  \, \times   \,   \# J_{n}    \qquad    \forall  \  n = 0,1, \ldots    \  .
%$$
\begin{eqnarray}\label{eq18sv3fc}
  \#  \JJ_{n+1 }  \  & = &   \  \# \II_{n+1}^{\, -}  \ - \  \sum_{t=1}^{d}  \,  {\rm  l.h.s. \ of \ }  (\ref{eq18svbbdd})   \nonumber \\[1ex]   \nonumber
 & \ge &   (R  -  2 \, \lceil R^{1-\alpha_{\min}}\rceil ) \, \#   \JJ_{n}   \  - \ 3 d\; R^{1-\epsilon } \,  \#  \JJ_{n}  \nonumber \\[2ex] \nonumber
& =  &  (R-5 dR^{1-\epsilon})  \   \#  \JJ_{n}  \qquad    \forall  \  n = 0,1, \ldots    \  .
\end{eqnarray}
The upshot, is that
$$
\# \JJ_{n}\ge (R-5d R^{1-\epsilon})^n \ge 1  \qquad    \forall  \  n = 0,1, \ldots    \
$$
and therefore
$$
\bigcap_{t= 1}^{d}  \bad_{c(t)} (i_t,j_t)    \cap \Theta    \ \supset   \
\bigcap_{n=1}^\infty
\bigcup_{J\in\JJ_n}J   \neq  \emptyset      \ .
$$

\noindent This establishes Corollary~\ref{tscc} in the case the pairs $ (i_t,j_t) $  satisfy \eqref{neq1sv}. In order to complete the proof in full, we need to deal with the pairs $(1,0)$ and $(0,1)$.

\subsection{Dealing with $(1,0)$ and $(0,1)$ }\label{endsec}

By definition,  $\bad(1,0)  =  \{ (x,y)\in \RR^2\,:\, x\in \bad \}$.   Thus, the condition that $\theta\in \bad$ imposed in Corollary~\ref{tscc} implies that
$$
\bad(1,0)\cap \L_{\theta}  \, =  \, \L_{\theta}\  \ .
$$
In other words, the pair $(1,0)$ has absolutely no effect when considering the intersection of any number of different $\bad(i,j)$ sets with $\L_{\theta}$ nor does it in anyway  effect  the
modified construction of~\S\ref{mbc}.
%
%In other words,  the pair $(1,0)$ has absolutely no effect when considering the intersection of any number of $\bad(i,j)$ sets with $\L_{\theta}$. For the same reason,  the pair $(1,0)$  does not in anyway  effect the above modified construction.
%

In order to deal with intersecting $\bad(0,1)$  with $\L_{\theta}$, we show that the pair $(0,1)$ can  be easily  integrated within the  modified construction.
To start with,  note that
$$
\bad (0,1)\cap \Theta=\{ (\theta,y)\in [0,1)^2\;:\; y\in \bad \}.
$$
With $c_1$ as in \S\ref{mbc}, let
\begin{equation}\label{defc10}
c:=\frac{c_1}{2R^2}   \; .
\end{equation}
For the sake of consistency with the previous section, for $n\ge 0$  let
$$
\cC(n):=\left\{p/q \in \QQ\;:\; R^{n-1}\le H(p/q)< R^n\right\}  \qquad   {\rm where  \ } \qquad
    H(p/q)  := q^2   \, .
$$
Furthermore, let $\Delta (p/q )$  be the interval centered at $(\theta,p/q)$ with length
$$
\left|\Delta (p/q)\right| := \frac{2c}{H(p/q)}   \, .
$$
%With reference to \S\ref{mbc},  suppose that $(i_t,j_t)$ is $(0,1)$ for some $t=1,\ldots,d$.  Since $\cC(n)=\emptyset $ for $n=0$ and $n=1$,   the following  analogue of \eqref{eq18svbbdd}   allows us to deal with  intersecting $\bad(0,1)$  with $\L_{\theta}$. \emph{For $R\ge 4$, we have that}
With reference to \S\ref{mbc},  suppose that $(i_t,j_t)$ is $(0,1)$
for some $t=1,\ldots,d$.  Since $\cC(n)=\emptyset $ for $n=0$, the
following  analogue of \eqref{eq18svbbdd}   allows us to deal with
the pair $(0,1)$ within the  modified construction. \emph{For $R\ge
4$, we have that}
\begin{equation}\label{defc10db}
\#\{I_{n+1}\in \II_{n+1}^{\, -}  \;  :  \;  \Delta(p/q)  \cap
I_{n+1}\neq\emptyset  {\rm \  \ for  \ some \   }   \ p/q \in \cC(n)  \} \;\le\;
 3 \;   \#  \JJ_{n}    \ .
\end{equation}
%This estimate is easily deduced from  the following two facts.
To establish this estimate we proceed as follows.  First note that in view of \eqref{defc10}, we have that
$$
\frac{|\Delta(p/q)|}{|I_{n+1}|}\le 1  \ .
$$
Thus, any single interval  $\Delta(p/q)$ removes at most three intervals $I_{n+1}$ from $\II_{n+1}$. Next, for any two rationals   $p_1/q_1,p_2/q_2\in \cC(n)$ we have that
$$
\left|\frac{p_1}{q_1}-\frac{p_2}{q_2}\right|  \, \ge \, \frac{1}{q_1q_2}  \, \ge  \,   R^{-n}  \,  > \,   c_1R^{-n}    \, .
$$
Thus, there is at most one  interval $\Delta(p/q)$  that can possibly intersect any given interval $J_n$ from $\JJ_n$.  This together with the previous fact establishes \eqref{defc10db}.

\section{Proof of Theorem \ref{thmlem3}}\label{triangsec}

\noindent Let $R \ge 2 $ be an integer. We start by showing that two parallel lines
from $\cC(n,l)$ can not intersect $J_{n-l}$. For any line $L(A,B,C)\in
\cC(n,l)$ we have that
$$
R^{\lambda l}  \, R^{-\frac{nj}{j+1}}   \;  \stackrel{\eqref{zeq1}}{<}    \; B^{-1}   \ .
$$
Thus,  if two parallel lines $L_1(A_1,B_1,C_1)$ and $
L_2(A_2,B_2,C_2)$ from $\cC(n,l)$ intersect $J_{n-l}$ 
we must have that
$$
  R^{2\lambda l}  \, R^{-\frac{2nj}{j+1}}
\, \leq   \,  \frac{1}{B_1B_2}   \, \leq \, 
 \left|\frac{C_1}{B_1}-\frac{C_2}{B_2}\right|
 \, \le \,    |J_{n-l}| \, =  \, c_1R^{-n+l}
\, .
$$
However, this is clearly false  since $c_1 < 1 < \lambda $ and $2j < j+1$.
%
%************************************************************************
%
%Let $R \ge 2 $ be an integer. We start by showing that two
%parallel lines from $\cC(n,l)$ can not intersect $J_{n-l}$. Suppose
%that two parallel lines $L_1(A_1,B_1,C_1),
%L_2(A_2,B_2,C_2)\in\cC(n,l)$ intersect $J_{n-l}$. Then the ratios
%$A_1/B_1$ and $A_2/B_2$ should coincide. Hence the distance between
%two parallel lines on $\Theta$ is bounded below by
%$$
%\left|\frac{C_1}{B_1}-\frac{C_2}{B_2}\right| \ge \frac{1}{B_1B_2}
%\stackrel{\eqref{zeq1}}{>} R^{2\lambda l}  \, R^{-\frac{2nj}{j+1}}.
%$$
%Thus, we must have that
%$$
%R^{2\lambda l}R^{-\frac{2nj}{j+1}}\le c_1R^{-n+l}  \, =  \,
%|J_{n-l}| \ .
%$$
%However, this is clearly false  since $\lambda>1$, $c_1 < 1 $ and $2j < j+1$.
%
%
%
%************************************************************************

Now suppose we have three lines $L_1, L_2$  and $L_3$  from $\cC(n,l)$
that intersect $J_{n-l}$ but do not intersect one another at a
single point.
  In view of the above discussion, the three lines $L_m=L(A_m,B_m,C_m)$ corresponding to $m =1,2 $ or $3$  can not be parallel to one another and therefore we must have three
distinct intersection  points:
$$
P_{12}=L_1\cap L_2  \, , \quad P_{13}=L_1\cap L_3   \quad\mbox{and
}   \quad  P_{23}=L_2\cap L_3   \ .
$$
Since $P_{12},P_{13},P_{23}$ are rational points in the plane, they
can be represented in the form
$$
P_{st}=\left(\frac{p_{st}}{q_{st}},\frac{r_{st}}{q_{st}}\right)
\qquad (1\le s< t\le 3)
$$
where
$$
\frac{p_{st}}{q_{st}}=\frac{B_sC_t-B_tC_s}{A_sB_t-A_tB_s}  \qquad {\rm  and }
\qquad \frac{r_{st}}{q_{st}}=\frac{A_sC_t-A_tC_s}{A_sB_t-A_tB_s} \ .
$$

\noindent In particular,  there exists an integer $k_{st} \neq 0 $ such that
\begin{equation*}\label{eq17sv}
k_{st}  q_{st}   =    A_s B_t  -  A_t B_s     \qquad {\rm and }  \qquad  k_{st}  p_{st} = B_sC_t- B_t C_s   \
\end{equation*}
and,  without loss of generality we can assume  that  $q_{st} > 0$.   On a slightly different note,  the  three intersection points  $Y_m:= L_m \cap J_{n-l} $  are obviously  distinct and it is easily verified that
$$
Y_m  =    \Big( \theta, \, \frac{A_m\theta +C_m}{B_m} \Big)   \qquad (1\le m \le 3)   \ .
$$

\noindent Let ${\rm T} (
P_{12}P_{13}P_{23}) $  denote the triangle subtended by the points
$P_{12},P_{23}$ and $P_{13}$.    Then twice the  area of the
triangle is equal to the absolute value of the determinant  $$ \det
\ := \ \left| \begin{array}{ccccc} 1 & p_{12}/q_{12}  &
r_{12}/q_{12} \\  \\ 1 & p_{13}/q_{13}  & r_{13}/q_{13} \\ \\ 1 &
p_{23}/q_{23}  & r_{23}/q_{23}
\end{array} \right|  \ . $$
It follows that,
\begin{equation}\label{eq17}
\area \,  {\rm T} (P_{12}P_{13}P_{23})\ge
\frac{1}{2q_{12}q_{13}q_{23}} \ .
\end{equation}

\noindent On the other hand, ${\rm T} (P_{12}P_{13}P_{23})$  is
covered by the union of triangles $ {\rm T} (Y_1Y_2P_{12})\cup {\rm
T} (Y_1Y_{3}P_{13})\cup {\rm T} (Y_2Y_{3}P_{23})$. Thus
$$\area \, {\rm T} (P_{12}P_{13}P_{23})\le
\area \, {\rm T} (Y_1Y_2P_{12})+\area \, {\rm T}
(Y_1Y_{3}P_{13})+\area \, {\rm T} (Y_2Y_{3}P_{23}).
$$
Without loss of generality,   assume that ${\rm T} ( Y_1Y_2P_{12})$
has the maximum area. Then
$$
\area \, {\rm T}( P_{12}P_{13}P_{23})\le 3\cdot \area \, {\rm T}(
Y_1Y_2P_{12})=\frac32
|Y_1-Y_2|\cdot\left|\theta-\frac{p_{12}}{q_{12}}\right|.
$$
Now observe that
\begin{eqnarray*}
c_1R^{-n+l}  \, = \,   | J_{n-l}  |    \ge |Y_1-Y_2|  & = &    \frac{|(A_1B_2-A_2B_1)\theta - (B_1C_2 - B_2C_1)|}{|B_1B_2|}  \\[2ex]
& = &   \frac{|k_{12} \, q_{12}\theta- k_{12} \, p_{12}|}{B_1B_2}  \\[2ex]
& \ge &  \frac{|q_{12}\theta-p_{12}|}{B_1B_2}   \ .
\end{eqnarray*}
Hence
$$
\area  \, {\rm T} (P_{12}P_{13}P_{23}) \ \le \  \frac32  \,
c_1^2R^{-2(n-l)}\frac{1}{q_{12}}B_1B_2  \, .
$$
Therefore, on combining with  \eqref{eq17} we have that
\begin{equation}\label{svteq2}
R^{2n}  \  \le \ 3 c_1^2 R^{2l}B_1B_2q_{13}q_{23}   \  .
\end{equation}
We now show that since $c_1$ satisfies  (\ref{eqc1}) and therefore
\begin{equation}\label{teq2}
4 c_1  R^{\lambda i} \, \le  \, 1  \; ,
\end{equation}
%\begin{equation}\label{teq2}
%3 c_1^2\cdot (3+R^{2\lambda i})<1  \ ,
%\end{equation}
the previous  inequality \eqref{svteq2} is in fact false.  As a consequence, the triangle ${\rm T} (P_{12}P_{13}P_{23})$ has zero
area and therefore can not exist. Thus,  if there are two or more lines from $\cC(n,l)$  that intersect
$J_{n-l}$ then they are forced to intersect one another at a single
point.

\vspace*{2ex}

On using  the fact that $q_{st}  \le    | A_s |B_t +  |A_t| B_s  $,
it follows that
\begin{eqnarray}\label{eq18}
 {\rm r.h.s \ of \  } \eqref{svteq2}  & \le  &     3  \, c_1^2R^{2l}B_1B_2 \, (|A_1|B_3+|A_3|B_1) \; (|A_2|B_3+|A_3|B_2) \nonumber \\[2ex]
      & = & 3  \, c_1^2R^{2l}B_1B_2 \, \Big(|A_1|B_3  |A_2| B_3 \, +  \,    |A_1|B_3 |A_3|B_2   \nonumber  \\
        & & \hspace*{20ex} + \,  |A_3|B_1   |A_2|B_3   \, + \,   |A_3|B_1  |A_3|B_2  \Big)    \ .
\end{eqnarray}

\noindent By making use of  \eqref{zeq2} and \eqref{zeq1}, it is
easily verified that
$$
|A_t|B_t  \ =\ |A_t|B_t^i\ B_t^j< R^{ni}\ R^{-\lambda jl}\
R^{n\frac{j^2}{j+1}}\ = \ R^{-\lambda jl} R^{\frac{n}{j+1}}  \, .
$$
In turn it follows that each of the first three  terms  associated
with  \eqref{eq18} is bounded above by
$$ 3 c_1^2R^{  2 l (1-(1+j)\lambda)  }R^{2n}  \ \stackrel{\lambda > 1 }{\le}   \  3c_1^2R^{2n}    \ . $$
Turning our attention to the fourth term, since $L_1, L_3\in \cC(n,l)$
we have via \eqref{zeq1} that $B_1\le R^\lambda B_3$. Therefore,
\begin{eqnarray*}
 3 c_1^2R^{2l}|A_3|^2B_1^2B_2^2  & \le  &    3c_1^2
R^{2l+2\lambda i}  \, (|A_3|B_3^i)^2  \, B_1^{2j}B_2^2 \\[1ex]
      & \le  & 3c_1^2
R^{2l(1-\lambda j-\lambda)+2\lambda i}R^{2n}  \\[1ex]
  & \stackrel{\lambda > 1 }{\le}    &  3c_1^2 R^{2\lambda i}R^{2n}  \ .
\end{eqnarray*}
On combining this with the estimate for the first three terms, we
have that
$$
 {\rm r.h.s \ of \  } \eqref{svteq2}  \ \le  \    R^{2n} (  9 c_1^2   +   3c_1^2 R^{2\lambda i}  ) \  <  \  R^{2n}   \,  12 c_1^2 R^{2\lambda i} \  \stackrel{ \eqref{teq2} }{< } \   R^{2n}  \ .
$$
Clearly this is not compatible with the  left hand
side of  \eqref{svteq2}   and
therefore we must have that   \eqref{svteq2} is  false.

\vspace*{4ex}

\noindent{\em Remark.  } It is evident from the proof that the  statement of Theorem \ref{thmlem3} is true for any fixed interval of length $|J_{n-l}| :=   c_1 \, R^{-(n-l)}$.

\vspace*{2ex}

\section{Preliminaries for Theorem \ref{maincount}  \label{preprelim}}

In this section, we make various observations and establish results  that are geared towards proving Theorem \ref{maincount}. Throughout,  $R \ge 2 $ is  an integer and for $n \in \NN $ and $\tau \in \RR^{>0}$ we let
$$
\J \,  =   \, \J(n,\tau)
$$
denote a generic interval contained within $\Theta $ of length $\tau  R^{-n}$. Note that the position of $\J$ within $\Theta $ is not specified. Also, for an integer $k \ge 0$, we let $\cC(n,k) $ denote the collection of lines from $\cC(n) $ with height satisfying the additional condition given by \eqref{eq1}; that is
$$
%\cC(n,\ast,k)
\cC(n,k)\, := \, \left\{L= L(A,B,C)  \in \cC(n) \, : \,   2^kR^{n-1}  \, \le \, H(A,B) \, < \, 2^{k+1}R^{n-1}   \right\}   \ .
$$
Trivially, for any $ l \geq 0$ we have that
$$
\cC(n,l,k)   \, \subset \, \cC(n,k)   \ .
$$
No confusion with the collection $\cC(n,l)$ introduced earlier in \S\ref{secdescrip} should arise. The point is that beyond  Theorem \ref{thmlem3},  the collection $\cC(n,l)$  plays  no further role in establishing Theorem \ref{maincount} and therefore will not be explicitly  mentioned.

\subsection{A general property}

The following is a general property concerning points in
the set  $\bad(i) $ and lines passing through a given rational
point in the plane.

\begin{lemma}\label{lem5}
Let  $\theta\in \bad(i) $ and $P:=(\frac{p}q,\frac{r}q)$ be a
rational  point such that
$$ |q\theta-p|  < c(\theta) \, q^{-i}  \ .$$
Then there exists a line $L=L(A,B,C)$  passing through P with
$|A|\le q^{i}$ and $0<B\le q^j$.
\end{lemma}

\noindent{\em Proof. \, }
Consider the set
$$ ap - br    \pmod{q}\quad\hbox{where}\quad 0\le a\le [q^i]\quad\hbox{and}\quad 0\le b\le [q^j]  \, .$$
The number of such pairs $(a,b)$ is
$$(q^i+1-\{q^{i}\})(q^j+1-\{q^j\})>q.$$
Therefore, by the `pigeon hole' principle,  there exist pairs
$(a_1,b_1)$ and $(a_2,b_2)$ such that
$$a_1p-b_1r\equiv a_2p-b_2r \pmod{q}.$$
Thus,  there is clearly a choice of integers $A,B,C$ with
$$Ap-Br+Cq=0\quad\hbox{where}\quad |A|\le q^i\quad\hbox{and}\quad 0\le B\le q^j.$$

\vspace{2ex}

It remains to show that we may choose $B>0$. This is where  the
Diophantine condition on $\theta$ comes into play. Suppose $B=0$.
Then $Ap+Cq=0$ and without  loss of generality,  we may assume
that $(A,C)=1$. Put $d:=(p,q)$ and define $q_*:=q/d $  and
$p_*:=p/d$. Then
$$Ap_*=-Cq_*\quad\hbox{and}\quad |A|=q_*  \  .$$
Hence   $q_*\le q^i $ and  $d\ge q^j   \ge {q_*}^{j/i}$. However
 $$d|q_*\theta   - p_*|=  |q\theta - p  |  \; <  \;  c(\theta)  \,  q^{-i}   \, . $$
Thus, it follows that
\begin{equation*}
|q_*\theta   -  p_*  |  \;< \; c(\theta) {q_*}^{-i}d^{-1-i}  \; \le \;
c(\theta)   q_*^{-1/i}  \, .
\end{equation*}

\noindent But this contradicts the hypothesis that $\theta\in
\bad(i)$ and so we must have that  $B>0$.
\newline\hspace*{\fill}$\boxtimes$

\vspace*{2ex}
%
% In the following subsections, we make various observations and establish results  that are specifically geared towards proving Theorem \ref{maincount}. Throughout,  $R \ge 2 $ is  an integer and for $n \in \NN $ and $\tau \in \RR^{>0}\!$,  we let
%$$
%\J \,  =   \, \J(n,\tau)
%$$
%denote a generic interval contained within $\Theta $ of length $\tau  R^{-n}$. Note that the position of $\J$ within $\Theta $ is not specified. Also, for an integer $k \ge 0$, we let $\cC(n,k) $ denote the collection of lines from $\cC(n) $ with height satisfying the additional condition given by \eqref{eq1}; that is
%$$
%%\cC(n,\ast,k)
%\cC(n,k)\, := \, \left\{L= L(A,B,C)  \in \cC(n) \, : \,   2^kR^{n-1}  \, \le \, H(A,B) \, < \, 2^{k+1}R^{n-1}   \right\}   \ .
%$$
%No confusion with the collection $\cC(n,l)$ introduced in \S\ref{secdescrip} should arise. The point is that beyond  Theorem \ref{thmlem3},  the collection $\cC(n,l)$  plays  no further role in establishing Theorem \ref{maincount} and therefore will not be  mentioned. Trivially, for any $ l \geq 0$ we have that
%$$
%\cC(n,l,k)   \, \subset \, \cC(n,k)   \ .
%$$

\subsection{Two non-parallel lines intersecting $\J(n,\tau)$}

Let  $P:=(\frac{p}{q},\frac{r}{q})$ be a rational point  in the plane and
consider two non-parallel lines
$$
\begin{array}{ll}
L_1:& A_1x-B_1y+C_1=0,\\[1ex]
L_2:& A_2x-B_2y+C_2=0  \ \
\end{array}
$$
that  intersect one another at  $P$. It follows that
% In view of  \eqref{eq17sv}
%It is easily verified that
$$
\frac{p}{q}=\frac{B_1C_2-B_2C_1}{A_1B_2-A_2B_1}  \qquad {\rm  and }
\qquad \frac{r}{q}=\frac{A_1C_2-A_2C_1}{A_1B_2-A_2B_1} \ .
$$
Thus, there exists an integer $t \neq 0$ such that
\begin{equation}\label{eq19}
A_1B_2-A_2B_1=tq   \qquad {\rm  and }  \qquad B_1 C_2 -  B_2C_1=tp  \,
.
\end{equation}
\noindent {\em Without loss of generality, we will assume that $ q >
0$.} In this section, we  investigate the situation in which  both lines pass through a generic interval  $\J   =    \J(n,\tau)$. Trivially, for this to happen  we must have that
$$
|\J|  \ \ge   \   |Y_1-Y_2|  \  =  \ \frac{|(A_1B_2-A_2B_1)\theta
-(B_1C_2-B_2C_1)|}{|B_1B_2|}
$$
where
$$
Y_m := L_m\cap \J  =  \Big( \theta \, , \;  \frac{A_m\theta +C_m}{B_m} \Big)  \quad \quad m=1,2  \ .
$$
This together with \eqref{eq19} implies that
\begin{equation}\label{eq4}
\frac{|q\theta-p|}{B_1B_2}  \ \le  \ \frac{|tq\theta-tp|}{B_1B_2}  \ \le  \  \tau \, R^{-n}   \, .
\end{equation}

\noindent In the case that the lines  $L_1 $ and $L_2  $ are from the collection  $\cC(n,k)$, this general estimate leads to the following statement.

\begin{lemma}\label{lem1}
%Let  $L_1,L_2 \in \cC(n,l,k)$ be  two lines that intersect
%at $P:=(\frac{p}{q},\frac{r}{q})$ and let  $\J   =    \J(n,\tau)$
%be a generic interval. Suppose
Let  $L_1,L_2 \in \cC(n,k)$ be  two lines that intersect at $P:=(\frac{p}{q},\frac{r}{q})$ and let  $\J   =    \J(n,\tau)$  be a generic interval. Suppose
$$
L_1\cap \J\neq \emptyset     \qquad { and }  \qquad L_2\cap
\J\neq\emptyset \ .
$$
Then
\begin{equation}\label{eq3}
|q\theta-p|  \ <   \ 2^i\tau  \, \frac{2^{k+1}}{R}  \, q^{-i}  \; .
\end{equation}
\end{lemma}

\noindent{\em Proof. \, }
With reference to  the lines $L_1= L(A_1,B_1,C_1)$  and $L_2=L(A_2,B_2,C_2)$, there is no loss of generality in assuming that $B_1\le B_2$.   With this mind, by \eqref{eq4} we have that

\begin{eqnarray} \label{yes}
|q\theta-p|  & < &  \tau R^{-n}B_1B_2 \ \nonumber \\[1ex]  & \stackrel{\eqref{eq1}}{\le}  & \tau  \,
R^{-n}  \, B_1  \, \left(2^{k+1}R^{n-1}\right)^{\frac{j}{j+1}}  \nonumber  \\[1ex] & = & \tau \, 2^{k+1} \, R^{-1}   \, B_1  \left(
2^{k+1}R^{n-1} \right)^{-\frac{1}{1+j}}   \; .
\end{eqnarray}
On the other hand, by \eqref{eq19} we have that
\begin{eqnarray} \label{yesyes}
q &\le&  |tq| \ =  \ |A_1B_2-A_2B_1|  \ \le  \  |A_1B_2|+|A_2B_1| \nonumber \\[2ex]
&\stackrel{\eqref{eq1}}{\le} &
\left(2^{k+1}  \, R^{n-1} \right)^{\frac{j}{j+1}} \left(\frac{2^{k+1}  \, R^{n-1}}{B_1}\right)^i   \ +  \  \left(2^{k+1}  \, R^{n-1} \right)^{\frac{j}{j+1}} \left(\frac{2^{k+1}  \, R^{n-1}}{B_2}\right)^i  \nonumber \\[2ex]
& =  &
 \left( B_1^{-i}  \, + \, B_2^{-i} \right)  \ \left(2^{k+1}  \, R^{n-1} \right)^{i+\frac{j}{j+1}} \nonumber \\[2ex]
&\le&  2  \, B_1^{-i}  \,
\left(2^{k+1}  \, R^{n-1} \right)^{\frac{1+ij}{1+j}}  \nonumber  \; .
\end{eqnarray}
Therefore
\begin{eqnarray*} \label{yesyesyes}
q^{-i}& \ge & 2^{-i}  \,  B_1^{i^2}  \,  \left(2^{k+1}  \, R^{n-1} \right)^{-\frac{i+i^2j}{1+j}} \nonumber  \\
& =  & 2^{-i}B_1
\ B_1^{-j(i+1)}  \,
\left(2^{k+1}  \, R^{n-1} \right)^{-\frac{i+i^2j}{1+j} }   \nonumber  \\
&\stackrel{\eqref{eq1}}{\ge} &2^{-i}  \, B_1  \,  \left(
2^{k+1}R^{n-1} \right)^{-\frac{j^2(i+1)}{1+j}-\frac{i+i^2j}{1+j}} \nonumber  \\
&=  &2^{-i}  \, B_1 \,  \left(
2^{k+1}R^{n-1} \right)^{-\frac{1}{1+j}}   \, .
\end{eqnarray*}

\noindent This estimate together with \eqref{yes} yields the desired statement.
\newline\hspace*{\fill}$\boxtimes$

\vspace*{4ex}

\noindent{\em Remark.  } It is evident from the proof that the  statement of Lemma \ref{lem1} is actually true for lines $L_1,L_2$ with height bounded above by $2^{k+1} R^{n-1}$.

\vspace*{2ex}

\subsection{The  figure $F$ } \label{secF}

In this section, we  give a geometric characterization of  lines from $\cC(n,l,k)$ that pass through a  given rational point and intersect a generic interval. Let  $L_1= L(A_1,B_1,C_1)$  and $L_2=L(A_2,B_2,C_2)$  be two lines from $\cC(n,l,k)$  that pass through $P:=(\frac{p}{q},\frac{r}{q})$ and intersect  $\J   =    \J(n,\tau)$.  Without loss of generality assume that  $B_1\le B_2$. Then, in view of  \eqref{eq4},  we have that
$$
\frac{|q\theta-p|}{B_1B_2}   \ \le  \  \tau R^{-n}  \ \stackrel{\eqref{eq1}}{<} \ \tau
\frac{2^{k+1}}{R}   \frac{1}{ H(A_2,B_2)}  \ .
$$
Thus
\begin{eqnarray} \label{meq1}
\frac{2^{k+1}\tau
}{R|q\theta-p|}  &  >  & \frac{H(A_2,B_2)}{B_1B_2}  \ = \ \frac{\max\{|A_2|^{1/i},B_2^{1/j}\}}{B_1}  \nonumber  \\[2ex]
& \ge &  \max \left\{\frac{|A_2|^{1/i}}{B_2} \, , \; B_2^{i/j}\right\}  \, .
\end{eqnarray}

\noindent Given a rational point $P$, the upshot  is that if two lines from $\cC(n,l,k)$  pass through $P$ and intersect $\J$,  then the point $(A,B) \in \ZZ^2$  associated with the   coordinates $A$ and $B$ of at least one of
the lines  lies inside the figure $F$ defined by

\begin{equation}\label{eq2}
|A|  \ < \ c_2^{i}B^i, \qquad 0 \, < \, B \, < \, c_2^{j/i} \qquad\mbox{with } \quad
c_2:=\frac{2^{k+1}\tau }{R|q\theta-p|}  \, .
\end{equation}

\vspace*{7ex}

\begin{center}
\begin{pspicture}(0,0)(6.5,-3)
\psset{linewidth=.2pt} \rput(3.4,-0.2){$c_2^{j/i}$}
\psline{->}(0,-2.5)(6.5,-2.5) \psline{->}(3,-3)(3,0.5)
\pscustom[fillstyle=vlines, hatchangle=45,hatchwidth=.2pt,
hatchsep=2pt]{
\pscurve(.5,-0.5)(1,-1.5)(2,-2.25)(3,-2.5)(4,-2.25)(5,-1.5)(5.5,-0.5)
\lineto(.5,-0.5)}
\rput(3.3,-2.7){O}\rput(6.4,-2.7){A}\rput(3.25,.5){B}
\psset{linestyle=dashed} \psline(5.5,-0.5)(5.5,-2.5)
\rput(5.3,-2.3){$c_2$}
\end{pspicture}
\end{center}

\begin{center}
The figure $F$ \vspace{10pt}
\end{center}

\vspace*{2ex}

\noindent Notice that the figure $F$ is independent of $l$ and
therefore the above discussion is actually  true for lines coming
from the larger collection $\cC(n,k)$. As a consequence, {\em apart
from one  possible   exception, all  lines   $L(A,B,C) \in \cC(n,k)$
passing through  $P$ and intersecting a generic interval~$\J$ will
have $A$ and $B$ coordinates corresponding to points $(A,B)$ lying
inside the figure~$F$}. Additionally, notice that the triple
$(A,B,C)$  associated with any  line $L$  passing through~$P$
belongs to the lattice

$$
\LLL=\LLL(P)   \ :=  \   \{ (A,B,C)\in
\ZZ^3 :   Ap-Br+Cq=0      \}  \ .  $$

\noindent We will actually be interested in the projection of $\LLL$
onto the $(A,B)$ plane within which the figure $F$ is embedded.  By an abuse of notation we will also refer
to this projection as $\LLL$.

\vspace*{2ex}

\noindent{\em Remark.  } Note that the figure $F$ is independent of
the actual position of the generic interval~$\J$.  However, it is
clearly  dependent  on the position of the rational point $P$.

\vspace*{2ex}

%\noindent{\em Remark 2.  } With reference to  $\cC(n,l,k)$, note that the figure
%$F$ is independent of $l$ and therefore the above discussion is true for the collection $\cC(n,k)$.
%In particular,  for an integer $k \ge 0$, apart from one  possible   exception, all
%lines   $L(A,B,C) \in  \cC(n,k)$  passing through  $P$ and intersecting
%$\J$ will give rise  to points $(A,B)$ lying within  $F$.

%\noindent{\em Remark 1.  } Note that the figure $F$ is independent of
%$n$ and the actual position of the generic interval $\J$.  It is
%dependent  only on the rational point $P$ and the parameters $\tau
%$, $R$ and $k$.

\vspace*{3ex}

Now assume that $L_1,L_2\in \cC(n,l,k)$  with  $l>0$. In this case we have that
\begin{eqnarray*}
B^{1/j} \ \stackrel{\eqref{zeq1}}{<} \ R^{-\frac{\lambda
l}{j}}R^{\frac{n}{j+1}} \ & \stackrel{\eqref{zeq3}}{<} & \ R^{1-\frac{\lambda
l(j+1)}{j}}  \, |A|^{1/i}  \\[2ex]
& \stackrel{\eqref{eq2}}{<} & \ R^{1-\frac{\lambda
l(j+1)}{j}}  \, c_2  \, B   \ .
\end{eqnarray*}
Therefore
\begin{equation}\label{zeq4}
0  \; <  \;  B \; <  \; c_3   \,  c_2^{j/i}  \  \qquad\mbox{with } \quad c_3:=R^{\frac{j}{i}-\frac{\lambda l(j+1)}{i}  }    .
\end{equation}

\noindent Note that $c_3 < 1 $ and  that
\begin{equation}\label{zeq5}
|A| \ \stackrel{\eqref{eq2}}{<} \ c_2^{i}B^i  \ <  \  c_3^i\cdot c_2  \, .
\end{equation}

\noindent The upshot  is that if two lines from $\cC(n,l>0,k)$  pass through $P$ and intersect $\J$,  then the point $(A,B) \in \ZZ^2$  associated with the   coordinates $A$ and $B$ of at least one of
the lines  lies inside the figure $F_l \subset F  $ defined by
\eqref{zeq4} and \eqref{zeq5}.

\subsection{Lines intersecting $\Delta(L_0)$      \label{rubbish} }
Let   $L_0=L(A_0,B_0,C_0)$  be  an arbitrary  line
passing through the rational point $P:=(\frac{p}{q},\frac{r}{q})$ and  intersecting $\Theta$. It is easily verified that the point $Y_0 :=  L_0 \cap \Theta $ has  $y$-coordinate
$$
 \frac{A_0\theta+C_0}{B_0} \ = \ \frac{A_0\frac{p}{q}+C_0}{B_0}+\frac{A_0}{B_0}\left(\theta-\frac{p}{q}\right)
 \ = \ \frac{r}{q}+\frac{A_0}{B_0}\left(\theta-\frac{p}{q}\right) \, .
$$
Now, assume there is another line $L=L(A,B,C)$ with
$$
H(A,B)  \ \ge  \  H(A_0,B_0)$$
passing through $P$ and intersecting $\Theta$.  Let
$$Y  =  Y(A,B,C)   \ :=  \  L \cap \Theta $$
and notice that
$$
Y\in \Delta(L_0) \ \Longleftrightarrow  \
|Y-Y_0| \ = \ \left|\frac{A}{B}-\frac{A_0}{B_0}\right|\left|\theta-\frac{p}{q}\right|  \ \le \  \frac{c}{H(A_0,B_0)} \ .
$$
In other words,
\begin{equation}\label{eq5}
Y\in \Delta(L_0)  \ \Longleftrightarrow \
\frac{A}{B}\in\left[\frac{A_0}{B_0}-\frac{c}{H(A_0,B_0)\left|\theta-\frac{p}{q}\right|},\frac{A_0}{B_0}+\frac{c}{H(A_0,B_0)\left|\theta-\frac{p}{q}\right|}\right].
\end{equation}
Geometrically, points  $(A,B)  \in \ZZ^2$ satisfying the right hand side of \eqref{eq5} form a cone $C(A_0,B_0)$ with
apex at origin. The upshot  is that all lines   $L=L(A,B,C)$ with  $A$ and $B$ coordinates
satisfying  $H(A,B)\ge H(A_0,B_0)$ and  $A/B  \in C(A_0,B_0)$, will  have
$Y(A,B,C)\in \Delta(L_0)$.

In addition, let $F$ be the figure associated with $P$, a generic interval $\J=\J(n,\tau)$ and the collection $\cC(n,k)$. Suppose that
\begin{equation}\label{neq6}
F\cap\LLL\subset C(A_0,B_0)  \qquad {\rm and }  \qquad  H(A,B)\ge
H(A_0,B_0)\quad \forall  \quad (A,B)\in F\cap\LLL.
\end{equation}
Then, in view of the discussion above,  any line $L=L(A,B,C)$
passing through $P$ such that   $(A,B)\in F\cap\LLL $  will have
$Y(A,B,C) \in \Delta(L_0)$. In particular, it follows via
\S\ref{secF} that  if we have two lines $L_1,L_2 \in \cC(n,k)$
passing through $P$ and intersecting  $\J$, then one of them has
coordinates corresponding to a point in $F\cap\LLL$ and therefore it
intersects $\J$  inside $\Delta(L_0)$. Thus, apart from one possible
exceptional line $L'$, all   lines $L=L(A,B,C) \in \cC(n,k)$ passing
through $P$ and intersecting  $\J$ will have the property that
$(A,B) \in F\cap\LLL $ and $Y(A,B,C) \in \Delta(L_0)$.
Note that for  $L'=L(A',B',C')$,  we have that $(A',B') \notin F\cap\LLL $ and therefore we can not  guarantee that  $H(A',B')\ge
H(A_0,B_0)$. Also, $L'$ may or may not intersects $\J$  inside $\Delta(L_0)$.

%WHAT DOES THIS MEAN FOR THICKENED L --  NAMELY IF $L_0$ IS FROM
%PREVIOUS LEVEL THEN   $\Delta (L)$  HAVE BEEN REMOVED IN
%CONSTRUCTION.

\subsection{The  key proposition}

Under the hypothesis of Lemma~\ref{lem1},  we know that there exists some  $\delta \in (0,1) $ such that
\begin{equation*}%\label{ceq1}
|q\theta-p|  \ =      \  \delta 2^i\tau \frac{2^{k+1}}{R} q^{-i}  \, .
\end{equation*}
Hence
\begin{equation}\label{meq9}
c_2  \stackrel{\eqref{eq2}}{=}  \frac{2^{k+1}\tau
}{R|q\theta-p|} \ = \ \delta^{-1}2^{-i}q^i \, .
\end{equation}

The following statement is at the heart of the proof of Theorem \ref{maincount}.

\begin{proposition}\label{lem2}
Let $P=(\frac{p}{q},\frac{r}{q})$ be a rational point  and $\J=\J(n,\tau)$ be a generic interval.    Let $\cC$ be the collection of lines $L=L(A,B,C)$ passing through $P$ with height $H(A,B) < R^n$. Let $\cC_k \subset \cC(n,k) $  denote the collection of lines passing through $P$ and intersecting $\J$.  Suppose that $\#\cC_k\ge 2 $,
$\tau \ge c R 2^{-k} $   and
\begin{equation}\label{lem2eq}
\delta\le c_4\left(\frac{cR}{2^k\tau }\right)^{2/j}  \qquad   {  where \ } \qquad   c_4
:=    4^{-2/j} \, 2^{-i}   \ .
\end{equation}
Then there exists a line $L_0\in \cC$ satisfying \eqref{neq6}. Furthermore,
apart from one possible exceptional line,  for all other  $L \in \cC_k$ we  have that
$(A,B) \in F\cap\LLL $ and $ L \cap \J \in \Delta(L_0)$.
%
%Then there exists a line $L_0\in \cC$ satisfying \eqref{neq6}. Furthermore,  let $\YY_k:= \{Y= L \cap \J :  L \in \cC_k\}$.  Then, apart from possibly one exceptional point $Y \in \YY_k$, we have that
%\begin{equation}\label{lem2eq2}
%\YY_k \subset \Delta(L_0)   \cup   Y  \,   .
%\end{equation}
\end{proposition}

\vspace*{2ex}

\noindent{\em Remark.  } We stress that the line $L_0$ of the proposition is completely  independent of the actual position of the generic interval $\J$ and therefore  the furthermore part of the proposition is also valid irrespective of the position of  $\J$.

\vspace*{2ex}

%\begin{proposition}\label{lem2}
%Let $P=(\frac{p}{q},\frac{r}{q})$ and $J$ be a some segment of
%length $\tau R^{-n}$ where $\tau \ge c_1$. Let $\mathbf{L}$ be the set
%consisting of lines of class at most $n$ going through $P$. Let
%$\mathbf{L}_k$ be the subset of $\mathbf{L}$ of all lines satisfying
%\eqref{eq1} and intersecting $J$. Let $\YY_k$ be the set of their
%intersection points with $J$.  Suppose that $\#\mathbf{L}_k\ge
%2 $ and $R > 1$. Then, if
%\begin{equation}\label{lem2eq}
%\|q\theta\|\le \frac{c}{2}
%\Bigl(\frac{r}{2^k}\frac{c}{\tau}\Bigr)^{\frac{1+i}{j}} q^{-i}
%\end{equation}
%there exists a line $L_0\in \mathbf{L}$  and possibly a point $Y \in
%\YY_k$ such that
%\begin{equation}\label{lem2eq2}
%\YY_k \subset \Delta(L_0)   \cup   Y  \,   .
%\end{equation}
%\end{proposition}
%
%Here the quantity  $\Bigl(\frac{r}{2^k}\frac{c}{\tau}\Bigr)$ is, for
%lines $L$ in \cC(n,l,k) (i.e satisfying (14))  $|J|/|\Delta(L)|$ so it
%counts the number of lines which could fit inside $ |J | $ and be disjoint.

\noindent{\em Proof.  } Notice that since $\#\cC_k\ge 2 $, there exists at least one line $L(A,B,C)  \in  \cC_k$ with   $A$ and $B$ coordinates corresponding to  $(A,B)$ lying within $F$  --  see  \S\ref{secF}.  Thus, there is at least one point  in $F\cap \LLL$ corresponding to a line with height bounded above by $R^n$.

\vspace*{2ex}

A consequence of
\S\ref{rubbish} is  that if there exists  a  line $L_0$ satisfying \eqref{neq6} then   the furthermore part of the statement of the proposition  is  automatically satisfied. In order to establish \eqref{neq6},
we consider the following two cases.

\vspace*{2ex}

\noindent {\bf Case A. \, }  Suppose there exists a point
$(A,B)\in F\cap \LLL$  such that
$$
B\le \sigma\cdot\delta\cdot q^j \qquad\mbox{where}\qquad
\sigma:=\left(2^{k+2+ij}\frac{\tau }{Rc}\right)^{1/j}   \, .
$$
Now let $(A'_0,B'_0)$ denote  such a point in $ F\cap \LLL$ with $B'_0$  minimal.
It follows that for all points $(A,B)\in F\cap\LLL$,
\begin{eqnarray*}
\left|\frac{A}{B}\right| \ & \stackrel{\eqref{eq2}}{<} &  \ \frac{c_2^i}{B^{1-i}}  \ \le \  \frac{(c_2 \,  B'_0)^i}{B'_0} \\[2ex]
& \le &  \frac{(\delta^{-1}2^{-i}q^i \   \sigma\delta
q^j)^i }{B'_0} \ =  \ \frac{2^{-i^2}\sigma^iq^i}{B'_0}    \
\end{eqnarray*}
and therefore
$$
\left|\frac{A}{B}-\frac{A'_0}{B'_0}\right|   \; <  \;  2 \ \frac{2^{-i^2}\sigma^iq^i}{B'_0}  \ .
$$

\noindent This together with \eqref{eq5}  implies  that if
\begin{equation}\label{eq12d}
\frac{c}{H(A'_0,B'_0)\left|\theta-\frac{p}{q}\right|}  \ \ge \
2 \, \frac{2^{-i^2}\sigma^iq^i}{B'_0}  \; ,
\end{equation}
then $F\cap\LLL\subset C(A'_0,B'_0)$. In other words,  the first condition of
\eqref{neq6} is satisfied.  Therefore, modulo \eqref{eq12d},   if the point $(A'_0,B'_0)$ has minimal height among all $(A,B)\in F\cap\LLL$ the second
condition of \eqref{neq6} is also  valid and  we are done. Suppose
this is not the case and let  $(A_0, B_0)$ denote the
minimal height point within $ F\cap\LLL$. Then $H(A_0,B_0)\le H(A'_0,B'_0)$ and so
\begin{eqnarray*}
\frac{c}{H(A_0,B_0)\left|\theta-\frac{p}{q}\right|} &\ge&  \frac{c}{H(A'_0,B'_0)\left|\theta-\frac{p}{q}\right|}     \\[2ex] &\stackrel{\eqref{eq12d}}{\ge} &
2 \, \frac{2^{-i^2}\sigma^iq^i}{B'_0} \ \ge \ \left|\frac{A}{B}-\frac{A_0}{B_0}\right|
\qquad \forall  \quad (A,B)\in F\cap\LLL  \, .
\end{eqnarray*}
Thus, by \eqref{eq5} we have that   $F\cap\LLL \subset C(A_0,B_0)$. The upshot is that if \eqref{eq12d} holds  then there exists a line from the collection $\cC$ satisfying \eqref{neq6}. We now establish~\eqref{eq12d}. Note that
\begin{eqnarray}
\eqref{eq12d} &  \Longleftrightarrow  &  \frac{c\cdot q^{1+i}}{
B'_0\max\{|A'_0|^{1/i},{B'_0}^{1/j}\}\cdot  2^i \tau \delta
\left(\frac{2^{k+1}}{R}\right)} \ \ge  \
\frac{2^{1-i^2}\sigma^iq^i}{B'_0}  \nonumber \\[3ex]
& \Longleftrightarrow & \left( \frac{cR}{2^{k+2+i-i^2}\tau \delta
\sigma^i}\right) q  \ \ge  \  \max\{|A'_0|^{1/i},{B'_0}^{1/j}\}  \, .  \label{eqstar}
\end{eqnarray}
Note that
$$
|A'_0|^{1/i} \ \stackrel{\eqref{eq2}}{<} \ c_2 B'_0   \, \le  \,   2^{-i}\sigma q  \ \quad {\rm and }   \ \quad  {B'_0}^{1/j}\le
\sigma^{1/j}  \delta^{1/j} q \, .
$$

\vspace*{2ex}

\begin{itemize}

\item Suppose that $|A'_0|^{1/i}>{B'_0}^{1/j}$.   Then
\begin{eqnarray*}
{\rm  r.h.s. \ of \ } \eqref{eqstar}  &  \Longleftarrow  & \frac{cR}{2^{k+2+i-i^2}\tau \delta \sigma^i} \ \ge  \
2^{-i}\sigma \\[2ex]
& \Longleftrightarrow &  \delta \ \le  \
\frac{cR}{2^{k+2-i^2}  \tau \sigma^{1+i}}   \ = \  c_4
\left(\frac{cR}{\tau 2^k }\right)^{\frac{2}{j}}  \ .
\end{eqnarray*}
%\frac{1}{R^{4/j}2^{j(1+i)+\frac{(1+i-i^2)(1+i)}{j}}}=
%\frac{1}{R^{4/j}2^\frac{(1+i)(1+j)}{j}}\mbox{ and}
%\end{equation}
This  is precisely  \eqref{lem2eq} and therefore verifies \eqref{eq12d} when $|A'_0|^{1/i}>{B'_0}^{1/j}$.

\vspace*{2ex}

\item Suppose that  $|A'_0|^{1/i}\le {B'_0}^{1/j}$. Then
\begin{eqnarray}  \label{eq9}
{\rm  r.h.s. \ of \ } \eqref{eqstar}  &  \Longleftarrow  &
 \frac{cR}{2^{k+2+i-i^2}\tau \delta \sigma^i}  \ \ge  \
\delta^{1/j}\sigma^{1/j} \nonumber \\[2ex]
&\Longleftrightarrow  &  \delta^{1+1/j}    \   \le \
\left(\frac{1}{2^{2+ij}}\right)^{1+\frac{ij+1}{j^2}}
\left(\frac{cR}{\tau2^k}\right)^{1+\frac{ij+1}{j^2}}  \nonumber \\[2ex]
&\Longleftrightarrow  &  \delta \ \le  \ c_4
\left(\frac{cR}{\tau 2^k}\right)^{1/j} \ .
\end{eqnarray}
By the hypothesis imposed on $\tau$, it follows that
\begin{equation}\label{eq9sv}
\frac{cR}{\tau 2^k} \le      1 \ .
\end{equation}
Therefore, in view of  \eqref{lem2eq} the
lower bound for $\delta $ given by \eqref{eq9} is valid.  In turn, this verifies  \eqref{eq12d} when $|A'_0|^{1/i}\le {B'_0}^{1/j}$.
\end{itemize}
\vspace*{2ex}

\noindent { \bf Case B. \, }
Suppose that for all points  $(A,B)$ within $F\cap\LLL$ we have that
$$
B  \ >  \  \sigma \delta  q^j \, .
$$
Then, in  view of  \eqref{eq2} it follows that

\begin{equation}\label{eq_pain}
\left|\frac{A}{B}\right| \ <  \  \frac{2^{-i^2}\sigma^i q^i}{\sigma\delta
q^j}   \ =  \ \frac{q^{i-j}}{2^{i^2}\sigma^j\delta}  \qquad \forall  \quad (A,B)\in F\cap\LLL  \, .
\end{equation}

\noindent By making use of  \eqref{svs},   \eqref{lem2eq}  and \eqref{eq9sv}, it is readily verified that
$$ |q\theta-p|  < c(\theta) \, q^{-i}    \ . $$
Thus, Lemma~\ref{lem5} is applicable and there exists a point
 $(A'_0,B'_0)  \in \LLL$ satisfying

$$
H(A'_0,B'_0)  \le     q^{1+ j } \, .
$$

\noindent As a consequence
\begin{equation}\label{veq1}
\frac{c}{H(A'_0,B'_0)\left|\theta-\frac{p}{q}\right|} \ \ge \
2\frac{q^{i-j}}{2^{i^2}\sigma^j\delta}  .
\end{equation}
Indeed,
\begin{eqnarray*}
 \eqref{veq1}  & \Longleftarrow &
\frac{cq^{1+i}}{2^{i}\tau\delta\left(\frac{2^{k+1}}{R}\right)
q^{1+j}}  \ \ge  \ \frac{q^{i-j}}{2^{i^2-1}\sigma^j\delta}  \nonumber \\[2ex]
& \Longleftrightarrow & \sigma^j  \ \ge  \ 2^{k+2+ij}\frac{\tau}{Rc}  \nonumber  \ .
\end{eqnarray*}
By the definition, the last inequality concerning $\sigma$ is valid and therefore so is \eqref{veq1}.  We now show that $F\cap\LLL \subset C(A'_0,B'_0)$.
 In view of \eqref{eq5}, this will be the case if
\begin{equation}\label{neq7}
\left|\frac{A}{B}-\frac{A'_0}{B'_0}\right|   \; \le  \;  \frac{c}{H(A'_0,B'_0)\left|\theta-\frac{p}{q}\right|} \qquad \forall  \quad (A,B)\in F\cap\LLL  \, .
\end{equation}

%To proceed, we consider two options.

\begin{itemize}
\item Suppose that  $(A'_0,B'_0)\in F\cap\LLL$.
%Then \eqref{neq7} is a consequence of \eqref{eq_pain} together  with \eqref{veq1}.
Then, clearly
$$\eqref{neq7}     \  \Longleftarrow \  \eqref{eq_pain}  \mbox{ \ and \ }  \eqref{veq1}   \, .
$$

\item Suppose that  $(A'_0,B'_0)\not\in F\cap\LLL$.  Then
\begin{eqnarray*}
 \eqref{neq7}  &  \Longleftarrow  &
\frac{c}{H(A'_0,B'_0)\left|\theta-\frac{p}{q}\right|} \ge \left|\frac{A}{B}\right|+ \left|\frac{A_0'}{B_0'}\right|
\\[2ex]  &  \Longleftarrow  &
\frac{c}{H(A'_0,B'_0)\left|\theta-\frac{p}{q}\right|}\ge
\left|\frac{A'_0}{B'_0}\right|+\frac{q^{i-j}}{2^{i^2}\sigma^j\delta}\\[2ex]
&\stackrel{\eqref{veq1}}{\Longleftarrow}&
\frac{c}{2H(A'_0,B'_0)\left|\theta-\frac{p}{q}\right|}\ge
\left|\frac{A'_0}{B'_0}\right|
\\[2ex]
&\Longleftarrow & \frac{cq^{1+i}}{2^i
\tau\delta\left(\frac{2^{k+1}}{R}\right) B'_0\, q}\ge
\frac{2q^i}{B'_0}\\[2ex]
& \Longleftrightarrow  &  \delta \le \frac{1}{4\cdot 2^i}  \
\frac{cR}{2^k  \tau}  \ .
\end{eqnarray*}
In view of  \eqref{lem2eq} and \eqref{eq9sv} this
lower bound for $\delta $ is valid and therefore so is \eqref{neq7}.
\end{itemize}

\noindent The upshot of the above is that $F\cap\LLL\subset C(A'_0,B'_0)$. In other words,  the first condition of \eqref{neq6} is satisfied. Therefore if the
pair $(A'_0,B'_0)$ has the minimal height among all $(A,B)\in
F\cap\LLL$ the second condition of \eqref{neq6} is also  valid and
we are done. Suppose this is not the case and let  $(A_0, B_0)\in
F\cap\LLL$ denote the minimal  height point within  $F\cap\LLL$. By assumption,
$$
H(A_0,B_0)< H(A'_0,B'_0)
$$
 and so
\begin{eqnarray*}
\frac{c}{H(A_0,B_0)\left|\theta-\frac{p}{q}\right|} &\ge&  \frac{c}{H(A'_0,B'_0)\left|\theta-\frac{p}{q}\right|}     \\[2ex] &\stackrel{\eqref{veq1}}{\ge} &
2 \, \frac{q^{i-j}}{2^{i^2}\sigma^j\delta} \ \
\stackrel{\eqref{eq_pain}}{\ge} \ \
\left|\frac{A}{B}-\frac{A_0}{B_0}\right| \qquad \forall  \quad
(A,B)\in F\cap\LLL  \, .
\end{eqnarray*}
Thus, by \eqref{eq5} we have that   $F\cap\LLL \subset C(A_0,B_0)$. The upshot is that  \eqref{neq7} holds  thus  there exists a line from the collection $\cC$ satisfying \eqref{neq6}.
\newline\hspace*{\fill}$\boxtimes$

\section{Proof of Theorem \ref{maincount} }\label{sec6}

Let $ l, k  \ge 0 $ and  $J_{n-l} \in \JJ_{n-l}$.
Let $\epsilon>0$ be sufficiently small  and $R=R(\epsilon)$ be sufficiently large.
In  view of the trimming process, Theorem \ref{maincount} will follow on showing  that no more than $R^{1-\epsilon} $   intervals  $I_{n+1}$ from $  \II_{n+1}$  can be removed by the intervals $ \Delta(L)$ arising from  lines  $L\in \cC(n,l,k)$ that intersect $J_{n-l}$. Let  $L_1,\ldots,L_M$ denote these lines of interest and let

$$
Y_m:=L_m\cap J_{n-l}  \ \qquad ( 1 \le m \le M)  \ .
$$

\noindent Indeed, then
$$ {\rm  l.h.s. \ of \ }  \eqref{countbadmain}  \  \le  \
\#\{I_{n+1}\in  \II_{n+1}^{\, - }  \;:\;
I_{n+1} \cap \Delta(L_m) \neq\emptyset    {\rm \  \ for  \ some \   }   \  1 \le m \le M  \}     \ .
$$
A consequence of   Theorem~\ref{thmlem3} is that the lines  $L_1,\ldots,L_M$ pass through a  single rational point $P=(\frac{p}{q},\frac{r}{q})$.  This is an absolutely crucial ingredient within the proof of Theorem \ref{maincount}.

\vspace*{2ex}

In view of \eqref{neq4} the number of intervals  $I_{n+1}  \in  \II_{n+1}$ that can be removed by any single  line $L_m$ is bounded above by
$$
\frac{2R^{n-\alpha}}{H(A,B)}+2  \  \stackrel{\eqref{eq1}}{\le}    \
K \, := \, \frac{2R^{1-\alpha}}{2^k}+2   \ .
$$

\noindent Notice that $K \ge 2 $ is independent of $l$. Motivated by the quantity $K$, we consider the following two cases.

\vspace*{2ex}

{\bf Case A.} \hspace*{2ex} Suppose that      $2^k<R^{1-\alpha}$.

\vspace*{1ex}

{\bf Case B.} \hspace*{2ex}  Suppose  that  $2^k\ge R^{1-\alpha}$.

\vspace*{2ex}

\noindent Then
$$
K  \ \le   \ \left\{\begin{array}{ll} \displaystyle
\frac{4R^{1-\alpha}}{2^k} \ &\mbox{ in  \ Case
A}\\[3ex]
\displaystyle 4  \ &\mbox{ in \ Case B}  \, .
\end{array}\right.
$$

\noindent Also, let
$$
\tilde{c}_1  \ :=  \ \left\{\begin{array}{ll} \displaystyle
\frac{4c_1R^{l+\epsilon-\alpha}}{2^k} \ &\mbox{ in \ Case
A}\\[3ex]
\displaystyle 4c_1R^{l+\epsilon-1}  \ &\mbox{ in  \  Case B}  \, .
\end{array}\right.
$$

\noindent We now subdivide the given  interval
$J_{n-l}$ into $d$  intervals $\tilde{I}_{nl}$ of equal length
%$c_1R^{l-n-1}\cdot K\cdot R^{\epsilon}\le \tilde{c}_1R^{-n}$
$c_1R^{l-n}\lceil R^{1-\epsilon} / K\rceil^{-1}$. It follows that
$$
d \ := \  \frac{|J_{n-l}|}{|\tilde{I}_{nl}|}  \ =  \ \left\lceil\frac{R^{1 -
\epsilon}}{K}\right\rceil
$$
and that
$$
|\tilde{I}_{nl}| \ :=  \ c_1R^{l-n}\lceil R^{1-\epsilon} / K\rceil^{-1} \ \le  \ \tilde{c}_1R^{-n} \, .
$$

\noindent By choosing $R $ sufficiently large and $\epsilon<\alpha$
so that
\begin{equation}\label{sv99}
R^{\alpha- \epsilon}  \ge  8 \, ,
\end{equation}
we can guarantee  that  \begin{equation}\label{sv999}  2 \  \le  \ d \  \le  \   \frac{2 \, R^{1 -
\epsilon}}{K}    \ .  \end{equation}

\noindent To proceed, we divide the $d$ intervals $\tilde{I}_{nl}$ into the following  two classes.

\vspace*{2ex}

{\bf Type 1.}  \hspace*{4ex}  Intervals $\tilde{I}_{nl}$ that intersect no more than one line among $L_1,\ldots,L_M$.

\vspace*{2ex}

{\bf Type 2.} \hspace*{4ex} Intervals $\tilde{I}_{nl}$ that intersect two or more  lines among $L_1,\ldots,L_M$.

\vspace*{3ex}

\subsection{Dealing with Type 1 intervals.}  Trivially, the number of  Type~1 intervals is bounded above by $d$. By definition, each Type~1 interval has no more than one line $L_m $ intersecting it. The total number of intervals  $I_{n+1} \in \II_{n+1} $  removed by a single line $L_m$  is bounded above by $K$.  Thus, for any strictly positive $\epsilon <  \alpha  $  and $R$ sufficiently large so that \eqref{sv99} is valid, the total number of intervals  $I_{n+1}  \in  \II_{n+1}$  removed by the lines $L_1, \ldots, L_M$  associated with Type 1 intervals is bounded above by

\begin{equation}\label{sveq18}
  d \,  K   \
\stackrel{\eqref{sv999}}{\le}   \   2\, R^{1-\epsilon}   \ .
\end{equation}

%Trivially, the number of  Type 1 intervals is bounded above by $d$. By definition, each Type 1 interval has no more than one line $L_m $ intersecting it. Therefore, the total number of intervals  $I_{n+1} \in \II_{n+1} $  removed from a single  Type 1 interval is bounded above by $K$. In particular, for any strictly positive $\epsilon <  \alpha  $  and $R$ sufficiently large so that \eqref{sv99} is valid, we have that
%
%\begin{eqnarray*}\label{sveq18}
%\#  \{I_{n+1}\in  \II_{n+1}  \!\!\! & : & \!\!\! I_{n+1} \cap \Delta(L_m)   \neq\emptyset {\rm \  \ and  \ } L_m  \cap  \tilde{I}_{nl}  \neq\emptyset      {\rm \  \ for  \ some \   }  \\ &  &   ~\hspace*{10ex} 1 \le m \le M {\rm \  \ and  \   Type \ 1 \ interval \ }  \tilde{I}_{nl}   \,  \}     \\[2ex]
%& \le &   d \,  K   \\[2ex]
%&\stackrel{\eqref{sv999}}{\le}  &   2\, R^{1-\epsilon}   \ .
%\end{eqnarray*}
%

\vspace*{3ex}

\subsection{Dealing with Type 2 intervals.}
%\noindent {\em Dealing with type (II) intervals.}
Consider an interval $\tilde{I}_{nl}$ of Type 2. By definition,
there are  at least two lines $L_s, L_t\in \cC(n,l,k)$ passing
through $P$ which intersect $\tilde{I}_{nl}$. With reference to
\S\ref{preprelim},  let $\J$ be a generic interval of length $
\tilde{c}_1R^{-n}$. Clearly $ |\J| $ is the same for $k$ and $l$
fixed and    $|\tilde{I}_{nl}| \le  |\J|  $. Also, in view of
\eqref{sv99}  we  have that  $  |\J|  < |J_{n-l}|  $. Thus, given an
interval $\tilde{I}_{nl}$ there exists a generic interval
$\J=\J(n,\tau)$ with $\tau := \tilde{c}_1$ such that $\tilde{I}_{nl}
\subset \J \subset   J_{n-l}$.
%******************************************
%
%DO WE NEED THAT $  |J|  < |J_{n-l}|  $ ? yes i think later when we
%show that no two points in ${F}\cap\LLL$ lie on a line passing
%through the origin.
%
%
%
%*****************************************************
By Lemma \ref{lem1}, there exists some  $\delta \in (0,1) $ such that
$$
|q\theta-p|  \; =  \; \delta2^i\tilde{c}_1q^{-i}   \,
\left(\frac{2^{k+1}}{R}\right).
$$
As a consequence of \S\ref{secF}, apart from one  possible   exception, all  lines   $L \in \cC(n,l,k)$ passing through  $P$ and intersecting $\J$ will have $A$ and $B$  coordinates corresponding to points $(A,B)$ lying inside the figure $F$   defined by \eqref{eq2} with
$c_2:=\delta^{-1}2^{-i}q^i$.  The upshot is that among the lines $L_1, \ldots , L_M $ passing through any $\tilde{I}_{nl}$ of Type 2,  all but possibly one line $L'$  will have coordinates corresponding to points in  $ F \cap \LLL$.  Moreover, if $l>0$ then  $F$ can be  replaced by the smaller figure $F_l$ defined by \eqref{zeq4} and \eqref{zeq5}.

%Consider all exceptional lines (i.e. the lines with the coordinates $(A,B)$ out of the figure $F$).

\vspace*{2ex}

\subsubsection{ Type 2 intervals with   $\delta$ small }
Suppose that
\begin{equation}\label{neq8}
\delta   \ \le  \  c_4 \left(\frac{cR}{2^k\tilde{c_1}}\right)^{2/j}   \qquad  {\rm  where \ } \qquad   c_4
:=    4^{-2/j} \, 2^{-i}   \ .
\end{equation}

\noindent With reference to the hypotheses of  Proposition
\ref{lem2}, the above guarantees \eqref{lem2eq}  and it is easily
verified that $  \tilde{c}_1   > c R 2^{-k} $ and that $\cC_k \ge 2$
since $\tilde{I}_{nl}  \subset \J$ is of Type 2.  Hence,
Proposition~\ref{lem2} implies the    existence of   a line $L_0 \in
\cC(n') $  with $ n' \le n$   passing through $P$ and
satisfying~\eqref{neq6}. Furthermore, among the lines $L_m$ from
$L_1, \ldots, L_M$ that intersect $\J$, all apart from  possibly one
exceptional line $L'$ will satisfy $ L_m \cap \J = Y_m \in
\Delta(L_0)$ and have coordinates corresponding to points $(A,B) \in
F\cap\LLL$. Note that $L_0$ is independent of the position of $\J$
and  therefore it is the same for each generic interval associated
with a Type 2 interval.  The point is that $P$ is fixed and all the
lines of interest pass though $P$.   However, in principle, the
possible exceptional line $L'$ may be different  for each  Type 2
interval.   Fortunately, it is easy to deal  with such lines.  There
are at most $d$ exceptional lines $L'$  --  one for each of the $d$
intervals $\tilde{I}_{nl}$. The number of intervals  $I_{n+1}  \in
\II_{n+1}$ that can be removed by any single  line $L'$ is bounded
above by $K$.  Thus, no more than $ d \, K  \leq 2\, R^{1-\epsilon}$
intervals $I_{n+1}$ are removed in total by the exceptional  lines
$L'$. Now consider those lines $L_m=L(A_m,B_m,C_m)$  among  $L_1,
\ldots, L_M$  that intersect some Type 2 interval  and are not
exceptional. It follows that
$$
Y_m  \in \Delta(L_0) \qquad {\rm and } \qquad  H(A_m,B_m ) \ge
H(A_0,B_0)  \ .
$$

\begin{itemize}

\item[$\circ$]   Suppose that  $L_0\in \cC(n')$ for some $n'<n$. Denote by $\Delta^+(L_0)$ the interval with the  same center as $\Delta(L_0)$ and length $|\Delta(L_0)| + 2\lceil R^{1+\alpha}\rceil |I_{n'+2}| $. It is readily verified that $\Delta(L_m)  \subset \Delta^+(L_0)$ for any non-exceptional line $L_m $. Now observe that the interval $\Delta(L_0)$ is removed (from the segment $\Theta$)  at level $n'$ of the basic construction; i.e. during the process of removing those `bad' intervals $I_{n'+1} $ from $ \II_{n'+1}^{\, -}$ that intersect some $\Delta(L)$ with  $L \in \cC(n')$.  The set $\Delta^+(L_0) \setminus \Delta(L_0)$ is removed (from the segment $\Theta$)  by the `trimming' process at level $n'+1$ of the basic construction.   In other words, the  interval $\Delta^+(L_0)$ has been totally removed  from $\Theta$ even before  we  consider the effect of lines  from $\cC(n)$ on the remaining  part of $\Theta$; i.e.  on intervals $I_{n+1}   \in  \II_{n+1}^{\, -}$.  In a nutshell, there are no intervals $I_{n+1}   \in  \II_{n+1}^{\, -}$ that lie in $\Delta^+(L_0)$ and therefore any non-exceptional line $L_m $  will have  absolutely  no `removal' effect.

\item[$\circ$] Suppose that  $L_0\in \cC(n)$.  Denote by $\Delta^+(L_0)$ the interval with the  same center as $\Delta(L_0)$ and length $2|\Delta(L_0)|$. It is readily verified that $\Delta(L_m)  \subset \Delta^+(L_0)$ for any non-exceptional line $L_m $. In view of \eqref{neq4}   the interval  $\Delta^+(L_0)$ can  remove no more than $4R^{1-\alpha}+2$ intervals  $I_{n+1}  \in  \II_{n+1}$.

\end{itemize}

\noindent The upshot when   $\delta$ satisfies
\eqref{neq8} is as follows.   For any strictly positive $\epsilon <  \alpha  $  and $R$ sufficiently large so that \eqref{sv99} is valid,
the total number of intervals  $I_{n+1}  \in  \II_{n+1}$  removed by the lines $L_1, \ldots, L_M$  associated with Type 2 intervals  is bounded above by
\begin{equation}\label{sveq19}
4R^{1-\alpha}+2 + K\cdot  d   \, =  \,
4R^{1-\alpha}+2+2R^{1-\epsilon} \, \le \, 6R^{1-\epsilon}+ 2   \, \le  \, 8
R^{1-\epsilon}  \ .
\end{equation}

 Naturally, we now proceed by dealing with the situation when  \eqref{neq8} is not satisfied.

\vspace*{2ex}

\subsubsection{ Type 2 intervals with   $\delta$ large }
Suppose that
\begin{equation}\label{neq8svsv}
\delta   \ >  \  c_4 \left(\frac{cR}{2^k\tilde{c_1}}\right)^{2/j}   \ .
\end{equation}

\noindent In Case A it follows   that
\begin{equation}\label{n40}
 \delta   \ >  \  c_4 4^{-2/j}  \, R^{-2(l+\epsilon)/j}
\end{equation}
and in Case B, using the fact that $2^k < R $ -- see \eqref{eq1vbvb},  it follows  that
\begin{equation}\label{n41}
\delta  \ > \  c_4 \left(\frac{R^{1-l-\alpha-\epsilon}}{4\cdot
2^k}\right)^{2/j}  \ >  \    c_4 4^{-2/j}
R^{-\frac{2(l+\alpha+\epsilon)}{j}}   \, .
\end{equation}

Recall that for the generic interval  $\J$ associated with a Type 2 interval $\tilde{I}_{nl}$,  there exists at most one exceptional line $L'$ among $L_1,\ldots,L_M$ that intersects $\J$ and has coordinates corresponding to a point not in  $F\cap\LLL$. We have already observed that no more than $ d \, K  \leq 2\, R^{1-\epsilon} $ intervals $I_{n+1}  \in  \II_{n+1}$ are removed in total by the $d$ possible exceptional lines $L'$. Indeed, the latter are exactly the same as in the $\delta$ small case  and therefore the  corresponding removed intervals  $I_{n+1} $ coincide.

We now consider those lines $L_m=L(A_m,B_m,C_m) $ among $L_1,\ldots,L_M$ with  $(A_m,B_m)\in F\cap\LLL$.
Suppose we have two such lines $L_m$ and $L_{m'}$ so that the points $(A_m,B_m)$ and $(A_{m'},B_{m'})$ lie on a line passing
through the lattice point $(0,0)$.  Clearly all points $(A,B)\in
F\cap \LLL$ on this line have the same ratio $A/B$. Thus the
lines $L_m$ and $L_{m'}$ are parallel. However, this is impossible since $L_m$ and $L_{m'}$ intersect at the rational point~$P$.  The upshot of
this is that the points $(A_m,B_m)$, $(A_{m'},B_{m'})$ and $(0,0)$   do not lie on the same line.    Recall that the lines $L_m$  of interest are from within  the collection $\cC(n,l,k)$.  To proceed we need to consider the $l=0$ and $ l > 0$ situations  separately.

\medskip

\noindent $\bullet$ {\rm Suppose that $l=0$}. Let $$M^*:=\#\{L_m \in\{L_1,\ldots,L_M\}\;:\; (A_m,B_m)\in F\cap\LLL\}.$$  Observe that the figure $F$ is convex.  In view of the discussion above, it then follows that  the lattice points in $F\cap \LLL$ together with the lattice point $(0,0)$ form
the vertices of  $(M^*-1)$ disjoint triangles lying within $F$. Since
the area of the fundamental domain of $\LLL$ is equal to $q$, the
area of each of these disjoint triangles is at least $q/2$  and
therefore the area of $F$ is at least $q/2\cdot (M^*-1)$.  Thus

$$
\frac{q}{2}(M^*-1)\le \area(F)< 2  c_2^{1+j/i}
\stackrel{\eqref{meq9}}{=}   \frac{q}{\delta^{1/i}}
$$
and therefore
$$
M^*   \, < \,    2 \, \delta^{-1/i}  +    1   \ .
$$
\begin{itemize}
\item[$\circ$] In Case A it follows via \eqref{n40} that
$$
M^*  \, <  \,     4^{\frac{4}{ij} +1   }    \;  R^{
\frac{2\epsilon}{ij} }  \, + \, 1 \ .
$$
Hence
$$
M^* K  < 20\cdot 4^{\frac{4}{ij}} \;
R^{1-\alpha+\frac{2\epsilon}{ij}}.
$$
Moreover, if  $$
\epsilon   \ \le \     \frac{\alpha \, ij }{ij + 2}  \,  ,
 $$ then we have that
\begin{equation}\label{ineq_mka}
M^* K   <  20 \cdot  4^{\frac{4}{ij}} \;   R^{1-\epsilon}      \ .
\end{equation}

\item[$\circ$]In Case B   it follows via \eqref{n41} that
$$
M^* <   4^{\frac{4}{ij} +1  }    \;  R^{ \frac{2(\alpha +
\epsilon)}{ij} }  \, + \, 1
$$
and thus the number of removed intervals is bounded by
$$
M^*K   <  20 \cdot  4^{\frac{4}{ij}} \;
R^{\frac{2(\alpha+\epsilon)}{ij}}  \, .
$$
It is readily verified that if
$$\epsilon   \ \le \     \frac{ij - 2 \alpha }{ij + 2}  \, , $$ then the upper bound for $M^*K$ given by \eqref{ineq_mka} is valid in Case B.
\end{itemize}
%Note that the value of $\epsilon$ doesn't depend on constants $R$ or
%$c_1$. It depends only on the pair $(i,j)$.

\vspace*{2ex}

\noindent $\bullet$ {\rm Suppose that $l>0$}.  Instead of working with the figure $F$ as in the $l=0$ situation,  we work with the `smaller' convex figure $ F_l \subset F$. Let
$$M^*:=\#\{L_m \in\{L_1,\ldots,L_M\}\;:\; (A_m,B_m)\in F_l\cap\LLL\}.$$  The same argument as in the $l=0$ situation  yields that
$$
\frac{q}{2}(M^* -1)  \ \le  \  \area(F_l)  \ <  \   2    c_3^{1+i}
c_2^{1+j/i}   \ =   \ R^{-\left(\frac{\lambda
l(j+1)}{j}-1\right)\cdot\frac{j(i+1)}{i}} \frac{q}{\delta^{1/i}}   \, .
$$

\begin{itemize}
\item[$\circ$]In Case A we have
$$
M^*   \ <  \   4^{\frac{4}{ij} +1}   \;
R^{\frac{2(l+\epsilon)}{ij}}\, R^{\frac{j(i+1)}{i}-\frac{\lambda
l(i+1)(j+1)}{i}} +  1   \, .
$$
Since $l > 0 $ and by definition $\lambda=3/j$, it follows that
\begin{equation}\label{ineq_lh0}
\frac{\lambda l(i+1)(j+1)}{i}-\frac{j(i+1)}{i}-\frac{2l}{ij}  \ > \ 0  \, .
\end{equation}
Thus
$$
M^*  \ < \  4^{\frac{4}{ij} +1 }  R^{\frac{2\epsilon}{ij}} +  1
$$
as in the  $l=0$ situation.  In turn, the upper bound for $M^*K$ given by \eqref{ineq_mka} is valid for  $l>0$.

\item[$\circ$]In Case B   we have
\begin{eqnarray*}
M^*& < & 4^{\frac{4}{ij} +1 } \,
R^{\frac{2(l+\alpha+\epsilon)}{ij}}\cdot
R^{\frac{j(i+1)}{i}-\frac{\lambda l(i+1)(j+1)}{i}} \; +  \; 1\\[2ex]
& \stackrel{\eqref{ineq_lh0}}{< } &   4^{\frac{4}{ij}+1 }   \,
R^{\frac{2(\alpha+\epsilon)}{ij}} \; + \; 1
\end{eqnarray*}
as in the  $l=0$ situation. In turn, the upper bound for $M^*K$ given by \eqref{ineq_mka} is valid in Case B for  $l>0$.

\end{itemize}

\noindent The upshot when   $\delta$ satisfies
\eqref{neq8svsv} is as follows.   For any strictly positive
$$
\epsilon   \ \le \     \frac{\alpha \, ij }{ij + 2}  \,  \stackrel{\eqref{eqc1sv}}{= }  \min \left\{  \frac{\alpha \, ij }{ij + 2}\, , \    \frac{ij - 2 \alpha }{ij + 2} \right\}
$$
and $R$ sufficiently large so that \eqref{sv99} is valid,
the total number of intervals  $I_{n+1}  \in  \II_{n+1}$  removed by the lines $L_1, \ldots, L_M$  associated with Type 2 intervals  is bounded above by
\begin{equation}\label{minsk}
K \, M^*    \ +  \  K\cdot  d   \, <  \, 20 \cdot  4^{\frac{4}{ij}} \;   R^{1-\epsilon}  \ +  \  2\,  R^{1-\epsilon}    \, <  \, 21 \cdot  4^{\frac{4}{ij}} \;   R^{1-\epsilon}   \ .
\end{equation}

\medskip

\subsection{The finale}
%To conclude,
On combining the upper bound estimates given by \eqref{sveq18}, \eqref{sveq19} and \eqref{minsk}, for any  strictly positive $\epsilon\le    \alpha  ij /(ij + 2)    $  and $R > R_0(\epsilon) $ sufficiently large, we have  that
\begin{eqnarray*}
 {\rm  l.h.s. \ of \ }  \eqref{countbadmain}  & \le&
\#\{I_{n+1}\in  \II_{n+1}^{\, - }  \;:\;
I_{n+1} \cap \Delta(L_m) \neq\emptyset    {\rm \  \ for  \ some \   }   \  1 \le m \le M  \}  \\[2ex]
& < & 2R^{1-\epsilon}+8R^{1-\epsilon}+  21 \cdot  4^{\frac{4}{ij}} \;   R^{1-\epsilon}     \ .
\end{eqnarray*}
This together with the fact that
$$
\alpha^2 \  <  \  \frac{\alpha \, ij }{ij + 2}
$$
completes the proof of Theorem  \ref{maincount}.

\newpage

\section{Proof of Theorem \ref{tsc}   \label{proofmainthm}  }

With reference to the statement of Theorem \ref{tsc}, since the set
under consideration is a subset of a line, we immediately obtain the
upper bound result that
\begin{equation}\label{uppbd}
 \dim \Big(\bigcap_{t=1}^{\infty} \bad(i_t,j_t)   \cap \L_{\theta}
\Big)  \le 1 \ .
\end{equation}

\noindent Thus, the proof of Theorem \ref{tsc} follows on
establishing  the following  complementary lower bound estimate.

\begin{theorem}\label{tsclb}
Let  $(i_t,j_t)$ be a countable number of pairs of real numbers
satisfying \eqref{neq1sv} and let $i := \sup\{i_t : t \in \NN \}$.
Suppose that
\eqref{liminfassump} is also satisfied.
Then, for any $\theta \in \bad(i)$ we have that
$$
\dim \Big(\bigcap_{t=1}^{\infty} \bad(i_t,j_t)   \cap \L_{\theta}
\Big)  \ge   1 \ .
$$
\end{theorem}

\vspace*{2ex}

\noindent {\em Remark. \, }  Strictly speaking, in order to deduce Theorem \ref{tsc}  we should replace  \eqref{neq1sv} by  \eqref{neq1} in the above statement of Theorem~\ref{tsclb}.   However,  given the arguments set out in  \S\ref{endsec},  the proof of Theorem~\ref{tsclb} as stated can easily be adapted to deal with the `missing' pairs $(1,0)$ and $(0,1)$.

\vspace*{2ex}

A general and classical method
for obtaining a lower bound for the Hausdorff dimension of an
arbitrary set is the following mass distribution  principle -- see
\cite[pg. 55]{falc}.

\begin{lemma}[Mass Distribution Principle]
 Let $ \mu $ be a probability measure supported on a subset $X$ of $ \RR$.
Suppose there are  positive constants $a, s $ and $l_0$ such that
\begin{equation}\label{mdp_eq1}
\mu  ( I ) \le \, a \;   |I|^s \; ,
\end{equation}
for any interval $I$ with length $|I|\le l_0$. Then,  $\dim X \ge s $.
\end{lemma}

\vspace*{2ex}

%Turning  to Theorem
%\ref{sc}, since the set under consideration
%is a subset of $\RR^2$, we immediately obtain the upper bound
%result that
%\begin{equation*}
% \dim \big(\bigcap_{t=1}^{\infty} \bad(i_t,j_t)
%\big)  \le 2 \ .
%\end{equation*}
%In view of (\ref{badi}) and  Theorem \ref{tsc},  the following
%general statement   enables us to establish the complementary
%lower bound estimate  -- see  \cite[pg. 99]{falc}.
%
%\begin{proposition}
%\label{vsections} Let $F$ be any subset of $\RR^2$, and let $E$ be
%a subset of the $x$-axis. If $\dim (F \cap \L_x)  \ge t $ for all
%$x \in E$, then $\dim  F  \ge  t + \dim E$.
%\end{proposition}
%
%\noindent Recall, that $\L_x$ denotes the line parallel to the
%$y$-axis passing through the point $(x,0)$ in the $(x,y)$-plane.

The overall strategy for establishing Theorem \ref{tsclb} is simply enough.  For each  $t\in \NN$,  let

\begin{equation}\label{666}
\alpha_t  \, :=  \,   \textstyle{\frac{1}{4}}  \, i_t \, j_t    \qquad\mbox{ and  } \qquad   \epsilon_0:  \, =   \,  \inf_{t\in \NN}\frac12 \alpha_t^2   \ .
\end{equation}

\noindent In view of  condition \eqref{liminfassump} imposed in the statement of the theorem,   we have that
 $\epsilon_0  $ is strictly positive. Then for any strictly positive $\epsilon<\epsilon_0$, we construct a `Cantor-type' subset $\KKK(\epsilon)$ of $ \bigcap_{t=1}^{\infty} \bad(i_t,j_t)   \cap \L_{\theta} $ and a probability measure $\mu$
supported on $\KKK(\epsilon)$ satisfying the condition that
\begin{equation}
\mu(I) \; \le \; a\, |I|^{1-\epsilon/2} \;   , \label{task}
\end{equation}
where the constant $a$ is absolute and $I \subset \Theta$ is an
arbitrary small interval.  Hence by construction and the mass distribution
principle we have that
$$
\dim \Big(\bigcap_{t=1}^{\infty} \bad(i_t,j_t)   \cap \L_{\theta}
\Big)  \; \ge \;   \dim \big( \KKK(\epsilon)  \big)   \; \ge  \;  1 - \epsilon/2  \ .
$$
Now suppose that $\dim \big(\bigcap_{t=1}^{\infty} \bad(i_t,j_t)
\cap \L_{\theta} \big)   < 1 $.  Then, $ \dim
\big(\bigcap_{t=1}^{\infty} \bad(i_t,j_t)   \cap \L_{\theta} \big) =
1 - \eta $ for some $\eta > 0$.  However, by choosing $\epsilon < 2
\eta $   we obtain a contradiction and thereby establish Theorem
\ref{tsclb}.

In view of the above outline, the whole strategy of our proof is
centred around the  construction of a `right type' of Cantor set
$\KKK(\epsilon)$ which  supports a measure $\mu$ with the desired
property.  It should  come as no surprise, that the first step involves  modifying  the basic construction to simultaneously incorporate any countable number of $(i,j)$ pairs satisfying \eqref{liminfassump} and  \eqref{neq1sv}.

\subsection{Modifying the basic construction for countable pairs}\label{sub_modconstr}

With reference to \S\ref{strategy},   for each $t \in \NN$ let us
 write $H_t(A,B)$  for $H(A,B)$, $\Delta_t(L)$ for $ \Delta(L)$   and $\cC_t(n)$ for $\cC(n)$. Furthermore, write $\JJ_n(t)$ for $\JJ_n$ and $\II^{-}_n(t)$ for $\II^{-}_n$.  With this in mind, let $R \ge  2$ be an integer. Choose $c_1(t)=c_1(R,t)$ sufficiently small
so that
%\begin{equation} \label{eqc1fc}
%c_1   \ \le \ \textstyle{\frac{1}{4}}   R^{-\frac{3i_t}{j_t}}   \qquad \forall \  \ 1 \le t \le d \,
%\end{equation}
\begin{equation}\label{ineq_c1t}
c_1(t)   \ \le \ \textstyle{\frac{1}{4}}   R^{-3i_t/j_t}   ,
\end{equation}
and
$$
c(t) \ := \   \frac{c_1(t)}{R^{1+\alpha_t}}
$$
satisfies \eqref{svs} with $\alpha_t  $ given by \eqref{666}.  With this choice of $ c_1(t)$,
the basic construction  of \S\ref{secdescrip} enables us to conclude
that $\bad(i_t,j_t)   \cap \L_{\theta} \neq \emptyset $ and in the process we establish the all important `counting' estimate  given by
\eqref{countbad}.  Namely,  let $ l \ge 0 $ and  $J_{n-l} \in
\JJ_{n-l}(t)$. Then, for any strictly positive $\epsilon <  \frac12
\alpha_t^2$ and $R
> R_0(\epsilon,t) $ sufficiently large we have that
\begin{eqnarray}  \label{999}
\#\{I_{n+1}\in \II_{n+1}^{\, -}(t)  \;  :  \; J_{n-l} \ \cap \   \Delta_t(L) \ \cap  \
I_{n+1}\neq\emptyset \   {\rm \  \ for  \ some \   }    \ L\in \cC_t(n,l)  \}  \ \le   \
 R^{1-\epsilon }    \ .
\end{eqnarray}
With $ \epsilon_0  $ given by \eqref{666}, this estimate is clearly valid for any strictly positive $\epsilon <  \epsilon_0 $.  The first step towards simultaneously dealing with the  countable number  of $(i_t,j_t)$ pairs  is to modify the basic construction in such a manner so  that corresponding version of  \eqref{999} remains intact.  The key is to start the construction with the   $(i_1,j_1)$ pair and then introduce at different levels within it the other pairs. Beyond this, the modifications are essentially at the `trimming' stage and in the manner in which  the collections $\JJ_{n}$  are defined.

% Suppose we have countably many pairs
%$(i_t,j_t)$, $t\in \NN$ with property \eqref{neq1}. We can do the
%main procedure for each pair independently i.e. we start with the
%segment $J_{0}(t)$ of length $c'_{1}(t)$ and on step $n$ of the
%procedure every segment is divided into $R'(t)$ pieces of equal
%length. Denote by $\JJ_n(t),\II_n(t)$ the sets analogous to
%$\JJ_n,\II_n$ defined in \ref{secdescrip}. Then, according to the
%procedure, lines $L(A,B,C)\in \cC(n,l,k)$ going through $J_{n-l}\in
%\JJ_{n-l}(t)$ can remove no more than $R'(t)^{1-\epsilon(i_t,j_t)}$
%subsegments from $\II_{n+1}(t)$.

Fix  some strictly positive $\epsilon  <   \epsilon_0$ and let $R$ be an arbitrary integer satisfying
\begin{equation}\label{eq20}
R   \, > \, R_0(\epsilon,1)  \       .
\end{equation}
Then, with
 $$c_1:=c_1(1)  \  $$
 we are able to carry out the basic construction  for the $(i_1,j_1)$ pair.   For each $t \ge  2$,  the associated basic construction for the $(i_t,j_t)$ pair is carried out with respect to a sufficiently large integer $R_t  > R_0(\epsilon,t) $ where $R_t$ is some power of $R$.   This enables us to embed  the  construction for each $t \ge 2$ within the construction for $t=1$.  More precisely, for $t\ge 1$ we let $$R_t:=R^{m_t}$$ where the integer $m_t$ satisfies
$$
m_1=1
$$
and for $ t \geq 2 $
$$  R^{m_t}\ge \max\{
R_0(\epsilon,t),   R^{1+m_{t-1}}   \} \ .
$$
Notice  that
\begin{equation}\label{ineq_mt}
m_t\ge t\ \mbox{ for }\ t\ge 2.
\end{equation}
Now  for each $t\ge 2$,  we fix an integer $k_t$ sufficiently large such that
$$
c_1(t):=c_1R^{-k_t}
$$
satisfies \eqref{ineq_c1t}  -- for consistency we let $k_1=0$.  Then for each $t \ge 1$, with this choice of $c_1(t)$ we are able to carry out the basic construction for the pair $(i_t,j_t)$. Moreover, for each integer $s\ge 0$ let
$$n_s(t):=k_t+sm_t   \, . $$
Then  intervals at  level  $s$ of the   construction for  $(i_t,j_t)$ can be described in terms of  intervals  at  level
$n_s(t)$ of the  construction for  $(i_1,j_1)$.  In
particular,  an interval of length $c_1(t)
R_t^{-s}$  at level $s$  for  $(i_t,j_t)$ corresponds to an
interval of length $c_1 R^{-n_s(t)}$  at level $n_s(t)$ for  $(i_1,j_1)$.

%In other words the basic construction for $(i_t,j_t)$ can be  incorporated within the basic construction for $(i_1,j_1)$ along the sequence given by $n_s(t)$.

We are now in the position to  modify the basic construction for the
pair $(i_1,j_1)$ so as to simultaneously incorporate  each
$(i_t,j_t)$ pair.  Let $c_1$ be as above. Define the collections
$\JJ_0:=\JJ_0(1)$ and $\JJ_1:=\JJ_1(1)$. Also Stage 1 of `the
induction' in which the collection  $\II_{n+1}$ is introduced
remains unchanged.  However the goal now is to remove those `bad'
intervals $I_{n+1}\in \II_{n+1}$ for which
\begin{equation}\label{badcountable}
I_{n+1}\cap \Delta_t(L)\neq\emptyset \ \mbox{ for some } t\in \NN\
\mbox{ and }  L\in
\cC_t\left(\textstyle{\left[\frac{n+1-k_t}{m_t}\right]-1}\right).
\end{equation}

%
%Regarding Stage~2, we trim the collection  $\II_{n+1}$ in the following manner.
%For any integer  $ t\ge 1$, if $n+1=n_{s+1}(t)$ for some $s$ we remove any interval $I_{n+1}$ from $\II_{n+1}$ that coincides with one of the $\lceil
%R_t^{1-\alpha_t}\rceil$ sub-intervals of length $|I_{n+1}|$ at either end of some $J_{n_s(t)}\in \JJ_{n_s(t)}$.  It follows that for  each such  $t$  the number of intervals $I_{n+1}$ from $\II_{n+1}$ that are removed by this `modified' trimming process is bounded above by
%$$
%\# \JJ_{n_s(t)}  \times 2 \, \lceil R_t^{1-\alpha_t}\rceil \, := \, \# \JJ_{n+1 -m_t}  \times 2 \, \lceil R^{m_t(1-\alpha_t)}\rceil \ ;
%$$
%i.e. the number removed by the basic  trimming  process associated with the pair $(i_t,j_t)$.  The intervals $I_{n+1}$ from $\II_{n+1}$ that survive the above trimming process  give rise to the trimmed collection $\II^-_{n+1}$. We define $\JJ_{n+1}$ to be the collection obtained by removing those `bad' intervals $I_{n+1}$ from $\II_{n+1}^-$ which satisfy \eqref{badcountable}.  In other words, for $n \ge  1$ we let
%\begin{equation}\label{svxz}
%\JJ_{n+1}:=\left\{I_{n+1}\in \II_{n+1}^-\;:\;
%\Delta_t(L)\cap I_{n+1}= \emptyset   \  \  \forall   \ t\in \NN    {\rm \  and \ } L\in C_t\left(\textstyle{\left[\frac{n+1-k_t}{m_t}\right]-1}\right)
%\right\}.
%\end{equation}
%Here, it is understood that the collection of lines
%$C_t(n)$ is the empty set whenever  $n$ is negative. Note that by construction,  the collection $\JJ_n$ is a sub-collection of $\JJ_s(t)$ whenever $n=n_s(t)$ for some $t\in\NN$.
%
%
%

Regarding Stage~2, we trim the collection  $\II_{n+1}$ in the
following manner. To begin with we remove from each $J_n\in\JJ_n$
the first $\lceil R^{1-\alpha_1}\rceil$ sub-intervals $I_{n+1}$ from
each end. In other words, we implement the basic trimming process
associated with the pair $(i_1,j_1)$.  Then for any integer  $ t\ge
2$, if $n+1=n_{s+1}(t)$ for some $s$ we incorporate the  basic
trimming process associated with the pair $(i_t,j_t)$. This involves
removing any interval $I_{n+1}$ that coincides with one of the
$\lceil R_t^{1-\alpha_t}\rceil$ sub-intervals of length $|I_{n+1}|$
at either end of some $J_{n_s(t)}\in \JJ_{n_s(t)}$.  It follows that
for  each such  $t$  the number of intervals $I_{n+1}$ from
$\II_{n+1}$ that are removed by this `modified' trimming process is
bounded above by
$$
\# \JJ_{n_s(t)}  \times 2 \, \lceil R_t^{1-\alpha_t}\rceil \, := \, \# \JJ_{n+1 -m_t}  \times 2 \, \lceil R^{m_t(1-\alpha_t)}\rceil \ ;
$$
i.e. the number removed by the basic  trimming  process associated
with the pair $(i_t,j_t)$.  Note that this bound is  valid for
$t=1$. The intervals $I_{n+1}$ from $\II_{n+1}$ that survive the
above trimming process  give rise to the trimmed collection
$\II^-_{n+1}$. We define $\JJ_{n+1}$ to be the collection obtained
by removing those `bad' intervals $I_{n+1}$ from $\II_{n+1}^-$ which
satisfy \eqref{badcountable}.  In other words, for $n \ge  1$ we let
\begin{equation}\label{svxz}
\JJ_{n+1}:=\left\{I_{n+1}\in \II_{n+1}^-\;:\;
\Delta_t(L)\cap I_{n+1}= \emptyset   \  \  \forall   \ t\in \NN    {\rm \  and \ } L\in C_t\left(\textstyle{\left[\frac{n+1-k_t}{m_t}\right]-1}\right)
\right\}.
\end{equation}
Here, it is understood that the collection of lines
$C_t(n)$ is the empty set whenever  $n$ is negative. Note that by construction,  the collection $\JJ_n$ is a sub-collection of $\JJ_s(t)$ whenever $n=n_s(t)$ for some $t\in\NN$.

%
%So there will be only finitely
%many nonempty classes
%$$C_t\left(\textstyle{\left[\frac{n+1-k_t}{m_t}\right]}\right)$$ for any natural
%$t$.

Apart from obvious notational modifications, Stages~3 and 4 remain
pretty much unchanged and gives rise to \eqref{999} for each $t\in
\NN$ with $R$ replaced by $R_t$. As consequence,  for any $ l\ge 0$
and $J_{n+1-(l+1)m_t}\in \JJ_{n+1-(l+1)m_t}$,  we have that
\begin{eqnarray}
\#\Big\{I_{n+1}\in\II^-_{n+1}  \!\!\! & : &   \!\!\!
J_{n+1-(l+1)m_t}\cap\Delta_t(L)\cap I_{n+1}\neq\emptyset  \nonumber     \\ & & \ \mbox{ for
some } \ L\in
C_t\left(\textstyle{\left[\frac{n+1-k_t}{m_t}\right]-1},l\right)\Big\}
  \ \le  \  R^{m_t(1-\epsilon)}   \, . \label{ineq_count}
\end{eqnarray}

\noindent To see this, let   $s+1 :=  \left[(n+1-k_t)/m_t \right] $.   Now if  $ s+1=  (n+1-k_t)/m_t $   then the statement is a direct consequence of  \eqref{999} with $ n  = s$ and $R$ replaced by $R_t$.  Here we use the fact that  $\II_{n+1}^-  \subseteq \II_{n+1}^-(t) $. Now suppose that  $  s+1 < (n+1-k_t)/m_t  $. Then $I_{n+1} $  is contained in some interval $J_{k_t + (s+1) m_t}$. By construction the latter does not intersect any  interval $\Delta_t(L)$  with $L \in \cC_t(s,l)$. Thus the set on the left hand side of \eqref{ineq_count}  is empty and the inequality is trivially satisfied.

\vspace*{2ex}

For fixed $\epsilon  < \epsilon_0 $ and any  $R$ satisfying \eqref{eq20}, the upshot of the modified basic construction is the existence of nested collections $\JJ_n$ of intervals $J_n$ given by \eqref{svxz} such that
\begin{equation}\label{obvious}
\KKK^*(\epsilon,R) \, :=   \, \bigcap_{n=0}^\infty \bigcup_{J\in \JJ_n} J  \ \subset \ \bigcap_{t=1}^\infty \bad(i_t,j_t)\cap \L_\theta  \, .
\end{equation}
Moreover, for $R$ sufficiently large  the counting estimate \eqref{ineq_count} can be used to deduce that
\begin{equation}\label{count_jn}
\#\JJ_n  \, \ge  \,  (R-R^{1-\epsilon/2})^n   \,
\end{equation}
-- see the remark following the proof of Lemma
\ref{sl_lem1} below.   Clearly, \eqref{count_jn} is more than
sufficient to conclude that $\KKK^*(\epsilon,R)$ is non-empty which
together with \eqref{obvious}  implies that
$$
\bigcap_{t=1}^\infty \bad(i_t,j_t)\cap \L_\theta \,  \neq  \emptyset \, .
$$
Recall, that as long as \eqref{liminfassump} is valid,  this enables us to establish the countable version of  Schmidt's conjecture.
However, counting alone is not enough to obtain the desired dimension result.  For this  we need to adapt the collections $\JJ_n$ arising from the modified construction.  The necessary `adaptation' will be the subject of the next section.

We end this section by investigation the distribution of intervals
within a given collection~$\JJ_n$.
%As a consequence  we shall establish  \eqref{count_jn}.
 Let $J_0$ be an arbitrary interval from $\JJ_0$ and define
$\TTT_0:=\{J_0\}$.  For $n \ge 1$, we construct the nested collections $ \TTT_n,  \TTT_{n-1},  \ldots, \TTT_1, \TTT_{0} $  as follows.  Take an arbitrary interval in $\TTT_{n-1} $ and subdivide it into $R$ closed intervals of equal length.  Choose
any $[2R^{1-\epsilon/2}]$ of  the $R$ sub-intervals and disregard the others.  Repeat this procedure for each interval in $\TTT_{n-1} $ and let $\TTT_n$ denote the collection of all chosen sub-intervals. Clearly,
$$
\# \TTT_n   \ =  \    \# \TTT_{n-1}  \times  [2R^{1-\epsilon/2}]  \, .
$$
Loosely speaking, the following result shows that the intervals $J_n $ from $\JJ_n$  are ubiquitous within each of  the intervals $J_0 \subset \Theta$ and thus within the whole of $\Theta$.
It is worth emphasizing  that both the collections $\JJ_n$ and $ \TTT_n$ are implicitly dependent on $R$.

 %Next, subdivide  $J_0$  into $R$ closed intervals of equal length and choose
%any $[2R^{1-\epsilon/2}]$ of them. Let $\TTT_1$   denote the
%resulting collection of intervals. We now repeat this procedure on each interval in  $\TTT_1$ and let $\TTT_2$ denote the resulting collection of intervals. We continue in this manner and construct the nested collections  $\TTT_3, \TTT_4, \ldots$.    Formally, we subdivide each interval from $\TTT_n$ into $R$ closed intervals of equal length and choose  any $[2R^{1-\epsilon/2}]$ of them. On $n$th step each
%segment from $\TTT_n$ is divided into $R$ equal pieces and some
%$[2R^{1-\epsilon/2}]$ of them are taken to $\TTT_{n+1}$.  Loosely speaking, the following result shows that the intervals $J_n $ from $\JJ_n$  are ubiquitous within the interval $J_0$ and thus within the whole of $\Theta$.
%

\begin{lemma}\label{sl_lem1} For $R$ sufficiently large,
\begin{equation}\label{sl_lem_eq}
\TTT_n\cap \JJ_n \neq \emptyset   \qquad \forall \quad  n = 0,1, \ldots \; .
\end{equation}
\end{lemma}

\noindent{\em Proof. \, }
For an integer $m \ge 0$, let $f(m)$ denote the cardinality of the set
$\TTT_m\cap \JJ_m$. Trivially,  $f(0)=1$ and the  lemma would follow on showing that
\begin{equation}\label{goodcount}
f(m)  \ \ge \  R^{1-\epsilon/2}  \, f(m-1)  \qquad \forall \quad  m \in \NN \, .
\end{equation}
This we now do via induction.
%It will easily imply the statement of the lemma. Indeed,
%$$
%\TTT\cap\bigcup_{J\in \JJ_n}J\ \supset\ \TTT\cap
%\bigcap_{n=1}^\infty\bigcup_{J\in \JJ_n}J\  =\
%\bigcap_{n=1}^\infty\bigcup_{J\in \TTT_n\cap \JJ_n} J.
%$$
%So we have the intersection of nested sets
%$$\bigcup_{J\in \KKK_n\cap \JJ_n} J$$ which in view of
%\eqref{ineq_hn} are nonempty.
To begin with, note that  $\# \JJ_1=\#\JJ_0 \times R $ and  so
\begin{equation*}%\label{ineq_f1}
f(1)\ =  \  [2R^{1-\epsilon/2}]  \ {>}   \  R^{1-\epsilon/2}   \, .
\end{equation*}
%Assume that \eqref{ineq_hn} is satisfied for $1,\ldots,n$. Then
%prove it for $n+1$. On $n$-th step each of the $f(n)$ intervals in
%$\TTT_n\cap \JJ_n$ gives rise to $[2R^{1-\epsilon/2}]$ intervals
%in $\TTT_{n+1}\cap \II_{n+1}$. In view of \eqref{ineq_count}, for any
%$t\in \NN, l\in\ZZ_{\ge 0}$ such that $n+1-(l+1)m_t\ge k_t$ there
%are at most
%$$f\big(n+1-(l+1)m_t\big)\cdot R^{m_t(1-\epsilon)}$$ segments from $\TTT_{n+1}\cap \II_{n+1}$ removed
%by lines from $C_t\big( [\frac{n+1-k_t}{m_t}]-1,l \big)$.
In other words,  \eqref{goodcount} is satisfied  for $m=1$. Now
assume that \eqref{goodcount} is valid for all $ 1 \le m \le n$. In
order to establish the statement for $m = n+1$, observe that   each
of the $f(n)$ intervals in $\TTT_n\cap \JJ_n$ gives rise to
$[2R^{1-\epsilon/2}]$ intervals in $\TTT_{n+1}\cap \II_{n+1}$. Now
consider some $t\in \NN$ and an integer $ l\ge 0$ such that
$n+1-(l+1)m_t\ge k_t$. Then in view of \eqref{ineq_count}, for any
interval $$J_{n+1-(l+1)m_t}\in
\JJ_{n+1-(l+1)m_t}\cap\TTT_{n+1-(l+1)m_t}$$ the number of intervals
from $\II^-_{n+1}$ removed by lines $L\in C_t\big(
[\frac{n+1-k_t}{m_t}]-1,l \big)$ is bounded above by
$R^{m_t(1-\epsilon)}$.  By the induction hypothesis,
$$
\# (\JJ_{n+1-(l+1)m_t}\cap\TTT_{n+1-(l+1)m_t})  \, =  \,  f(n+1-(l+1)m_t)   \ .
$$
Thus the total number of intervals from $\TTT_{n+1}\cap \II^-_{n+1}$ removed
by lines from $C_t\big( [\frac{n+1-k_t}{m_t}]-1,l \big)$ is bounded above  by
$$ R^{m_t(1-\epsilon)}  \  f\big(n+1-(l+1)m_t\big)    \ .$$

\noindent Furthermore, the number of intervals from $\TTT_{n+1}\cap \II_{n+1}$ removed by the modified trimming process associated with the pair $(i_t,j_t)$ is bounded above by

$$
2 \, \lceil R^{m_t(1-\alpha_t)}\rceil    \ f(n+1-m_t)  \ \le   \  2 \, R^{m_t(1-\epsilon)}  \ f(n+1-m_t).
$$
Here we have made  use of  the fact that $R^{m_t}  > R_0(\epsilon,t)$ and so  $ \lceil R^{m_t(1-\alpha_t)}\rceil   \leq R^{m_t(1-\epsilon)} $.

On combining the above estimates for intervals removed by `lines' and those removed by `trimming', it follows that
\begin{eqnarray*}
f(n+1)  &  \ge     &   [2R^{1-\epsilon/2}] \ \ f(n)  \\[2ex] & &  \hspace*{1cm} -  \ \ \sum_{t=1}^\infty
R^{m_t(1-\epsilon)}\sum_{l=1}^\infty f(n+1-l m_t)     \ - \   2\sum_{t=1}^\infty R^{m_t(1-\epsilon)}f(n+1-m_t) \  .
\end{eqnarray*}
Here, it is understood  that $f(k)=0$  whenever  $k$ is negative.
Then,  in view of our induction  hypothesis,  we have that
\begin{eqnarray*}
f(n+1) & \ge  &  [2R^{1-\epsilon/2}] \ \ f(n)   \\[2ex]
& & \hspace*{1cm} -  \ \ \sum_{t=1}^\infty R^{m_t(1-\epsilon)} \
f(n) \ (R^{-1+\epsilon/2})^{m_t-1}\left(2+\sum_{l=0}^\infty
(R^{-1+\epsilon/2})^{lm_t}\right)   \\[2ex]
& \ge & f(n)\left([2R^{1-\epsilon/2}] \ - \ R^{1-\epsilon} \ C(R)  \
-  \ \sum_{t=2}^\infty R^{1-\frac{m_t}{2}\epsilon  -
\frac{\epsilon}{2}} \ C(R) \right)
\end{eqnarray*}
where
$$
C(R):=2+\sum_{k=0}^\infty (R^{-1+\epsilon/2})^k.
$$
In addition to $R$ satisfying \eqref{eq20}  we assume that $R$ is sufficiently large so that
\begin{equation}\label{moreonR}
C(R) \ < \ 4 \, ,    \hspace*{4ex}   [2R^{1-\epsilon/2}] \ \ge  \ {\textstyle{\frac{5}{3}}} R^{1-\epsilon/2}\ \quad \mbox{ and }\
\quad \sum_{k=1}^\infty R^{-\frac{k}{2}\epsilon}  \, <  \, {\textstyle{\frac{1}{6}}}
\, .
\end{equation}
Then, by making use of \eqref{ineq_mt}  it follows that
$$
f(n+1)  \ \ge  \ R^{1-\epsilon/2}f(n) \left({\textstyle{\frac{5}{3}}} \ -
\ 4\sum_{t=1}^\infty R^{-\frac{t}{2}\epsilon}\right)  \ \ge \
R^{1-\epsilon/2}f(n)  \, .
$$
This completes the proof of the lemma.
\newline\hspace*{\fill}$\boxtimes$

\vspace*{2ex}

\noindent{\em Remark.  } For any $R$ satisfying \eqref{eq20} and \eqref{moreonR}, a straightforward consequence of \eqref{goodcount} is that
$$
\# \JJ_n     \, \ge \,    f(n)  \times \# \JJ_0  \   \ \ge \  R^{1-\epsilon/2}  \, f(n-1)  \times \# \JJ_0 \ \ge \  (R^{1-\epsilon/2})^n  > 1    \, .
$$
This is sufficient to show  that  $\KKK^*(\epsilon,R)$ is non-empty  and in turn enables us to establish the countable version of  Schmidt's conjecture.
However,  the proof of the lemma can be naturally modified adapted to deduce the  stronger counting estimate given by \eqref{count_jn} -- essentially replace $f(m) $ by $ \# \JJ_m  $   and $  [2R^{1-\epsilon/2}] $ by $R$.

\vspace*{2ex}

\subsection{The set $\KKK(\epsilon)$ and the measure $\mu$ }

Fix some strictly positive $\epsilon < \epsilon_0$ and an integer
$R$ satisfying \eqref{eq20} and \eqref{moreonR}.   The modified
construction of the previous section enables us to conclude that the
set $\KKK^*(\epsilon) := \KKK^*(\epsilon,R)$ defined by
\eqref{obvious} is  non-empty and in turn implies the weaker
non-empty analogue of Theorem~\ref{tsclb}.    To obtain the desired
dimension statement we  construct  a regular `Cantor-type' subset
$\KKK(\epsilon)$ of $ \KKK^*(\epsilon)$ and  a measure $\mu$
satisfying \eqref{task}.    The key  is to refine the  collections
$\JJ_n$ arising from the modified construction in such a manner that
the refined nested  collections $\MM_n  \subseteq \JJ_n  $ are
non-empty and  satisfy the following property. For any integer  $n
\ge 0 $ and $J_n\in \MM_n$
\begin{equation*}
\# \{J_{n+1} \in \MM_{n+1} \;:\; J_{n+1} \subset
J_{n}\}   \ \ge \   R-2R^{1-\epsilon/2}  \ .
\end{equation*}
%We have constructed a sequence of nested collections $\JJ_n$ however
%the resulting set $$\KKK=\bigcap_{n=1}^\infty \bigcup_{J\in\JJ_n}J$$
%is not a standard Cantor-type set. More precisely, we do not
%guarantee that each surviving interval $J_n\in \JJ_n$ contains at
%least one surviving interval $J_m\in\JJ_m$ for each $m>n$, so the
%intersection $\KKK\cap J_n$ is not necessarily nonempty. To be able
%to get the Hausdorff dimension result we will modify the collections
%$\JJ_n$.
%
%
%We will construct a subcollections $\MM_n\subset \JJ_n$
%with the following condition: for any $J_n\in \MM_n$ there are at
%least $R-2R^{1-\epsilon/2}$ segments $J_{n+1}\in\MM_{n+1}$ contained
%in $J_n$.
%
Suppose for the remaining part of this section  the desired  collections $\MM_n $  exist and let
$$
\KKK(\epsilon) \ := \ \bigcap_{n=0}^\infty \bigcup_{J\in \MM_n} J    \ .
$$
We  now construct a  probability measure $\mu$
supported on $\KKK(\epsilon)$ in the standard manner. For
any $J_n \in \MM_n$, we attach a weight $\mu(J_n)$
defined recursively as follows.

\vspace*{2ex}

\noindent For  $n=0$, $$ \mu(J_0)\ := \
\frac{1}{\# \MM_0}  \ $$ and
for $n\ge 1$,
\begin{equation}\label{beq4}
\mu(J_n)  \, :=  \, \frac{\mu(J_{n-1})}{\# \{J\in \MM_n\;:\; J\subset
J_{n-1}\}} \
\end{equation}

\vspace*{2ex}

\noindent where $J_{n-1} \in \MM_{n-1}$   is the unique interval such that $J_n\subset J_{n-1}$.
 This procedure thus defines inductively a mass on any interval
appearing in the construction of $\KKK(\epsilon) $. In fact a lot
more is true  ---   $\mu$  can be further extended to all Borel
subsets $F$ of $\RR $  to determine $\mu(F) $  so that $\mu$
constructed as above actually defines  a measure supported on
$\KKK(\epsilon) $.  We now state this formally.

\begin{itemize}
\item[]  {\em Fact.} The probability measure $\mu$ constructed above is
supported on $\KKK(\epsilon) $ and for any Borel set $F$
\[
\mu(F):= \mu(F \cap \KKK(\epsilon) )  \; = \;
\inf\;\sum_{J\in  \cJ }\mu(J)  \ .
\]
The infimum is over all coverings $\cJ$ of $F \cap
\KKK(\epsilon)$ by intervals  $J\in \{   \MM_n: n=0,1, \ldots \}$.
\end{itemize}

\vspace*{1ex}

\noindent For further details see  \cite[Prop. 1.7]{falc}.
It remains to show that $\mu \, $ satisfies \eqref{task}.  Firstly, notice that for any interval $J_n \in \MM_n$  we have that
\begin{eqnarray*}
\mu(J_n)  & \le  &
\left(R  \, \big(1-2R^{-\epsilon/2}\big)\right)^{-1}   \ \mu(J_{n-1}) \\[1ex]
& \le  &
\left(R\big(1-2R^{-\epsilon/2}\big)\right)^{-n}  \, .
\end{eqnarray*}
Next, let $d_n$ denote the length of a generic interval $J_n  \in \MM_n$  and consider an arbitrary interval $I\subset \Theta $  with length $|I| < d_0$. Then there exists a non-negative  integer $ n $ such that
\begin{equation} \label{fine}
d_{n+1}   \, \le  \, |I|  \,  <  \, d_n  \; .
\end{equation}
It follows that
\begin{eqnarray*}
\mu(I)  & \le  &  \sum_{\begin{array}{c}\scriptstyle J_{n+1}\in \MM_{n+1}\\
\scriptstyle J_{n+1}\cap I\neq\emptyset\end{array}} \mu(J_{n+1})   \\[2ex]  & \le  & \left\lceil \frac{|I|}{d_{n+1}}\right\rceil
 \big(R  \, (1-R^{-\epsilon/2})\big)^{-n-1}   \\[2ex]
& \le  & 2  \, \frac{|I|}{c_1 R^{-n-1}}  \  R^{-n-1}  \
\big(1-2R^{-\epsilon/2}\big)^{-n-1}  \\[2ex]
& \stackrel{\eqref{fine}}<  & 2 \,
c_1^{\epsilon/2-1}R^{\epsilon/2}\left(R^{\epsilon/2}\big(1-2R^{-\epsilon/2}\big)\right)^{-n-1}
|I|^{1-\epsilon/2}   \\[2ex]
&  \stackrel{\eqref{moreonR}}{\le}  & 2 \,
c_1^{\epsilon/2-1}R^{\epsilon/2}|I|^{1-\epsilon/2}   \ .
\end{eqnarray*}
Thus  \eqref{task} follows with
$a=2c_1^{\epsilon/2-1}R^{\epsilon/2}$ and this completes the proof
of Theorem \ref{tsclb} modulo the existence of the collection
$\MM_n$.

\subsection{Constructing the collection  $\MM_n$}

For any integer $n \geq 0$, the goal of  this section is to  construct the  desired nested collection $\MM_n  \subseteq \JJ_n$ alluded  to in the previous section.    This will involve constructing
auxiliary collections $\MM_{n,m}$ and $\RRR_{n,m}$ for  integers $n,m$  satisfying $0\le n \le m$. For a fixed $m$,  let
$$
\JJ_0  \, ,  \ \JJ_1 \, ,  \  \ldots, \ \JJ_m   \;
$$
be the collections arising from  the  modified construction of \S\ref{sub_modconstr}. We will require $\MM_{n,m}$ to satisfy the following
conditions.
\begin{itemize}
\item[\bf C1.] For any $0\le n\le m$, we have that $\MM_{n,m}\subseteq \JJ_n$.
\item[\bf C2.] For any $0\le n< m$, the  collections $\MM_{n,m}$ are nested; that is
    $$\bigcup_{J\in \MM_{n+1,m}}J  \qquad \subset \quad \bigcup_{J\in \MM_{n,m}}J.$$
\item[\bf C3.] For any $0\le n<m$ and $ J_n\in \MM_{n,m}$, we have that  there are at least  $R-[2R^{1-\epsilon/2}]$  intervals $J_{n+1}\in \MM_{n+1,m}$ contained within $J_n$; that is
\begin{equation*}
\# \{J_{n+1} \in \MM_{n+1,m} \;:\; J_{n+1} \subset J_{n}\}   \ \ge \
R-[2R^{1-\epsilon/2}]  \ .
\end{equation*}
\end{itemize}

\noindent  In addition, define
$\RRR_{0,0} := \emptyset $ and for $m \ge 1$
%$$
%\RRR_{0,j}:=\JJ_0\backslash \MM_{0,j}
%$$
\begin{equation}\label{def_rnn}
\RRR_{m,m}:=\left\{I_{m}\in \II_m\backslash \JJ_m\;:\; I_m\subset
J_{m-1}\mbox{ for some }J_{m-1}\in \MM_{m-1,m-1}\right\}   \ .
\end{equation}
Furthermore, for $ 0\le n< m $ define
\begin{equation}\label{def_rnm}
\RRR_{n,m}:=\RRR_{n,m-1}\cup \{J_n\in \MM_{n,m-1}\;:\;
\#\{J_{n+1}\in \RRR_{n+1,m}\;:\; J_{n+1}\subset J_n\}\ge
[2R^{1-\epsilon/2}]\ \} \ .
\end{equation}
Loosely speaking and  with reference to condition (C3),  the collections $\RRR_{n,m}$ are the `dumping
ground'  for those intervals $J_n\in\MM_{n,m-1}$ which do not contain
enough  sub-intervals $J_{n+1}$. Note that for $m$ fixed,
the collections $\RRR_{n,m}$ are defined  in descending order with respect to~$n$. In other words,  we  start with
$\RRR_{m,m}$  and finish with $\RRR_{0,m}$.

\vspace*{1ex}

The construction is as follows.

\vspace*{1ex}

\noindent{\em Stage 1.} Let $\MM_{0,0}:=\JJ_0  $  and   $ \RRR_{0,0}:=\emptyset$.

\vspace*{1ex}

\noindent {\em Stage 2.} Let $  0 \le t \le n $. Suppose  we have constructed the desired collections
$$\MM_{0,t}\subseteq \JJ_0, \   \MM_{1,t}\subseteq \JJ_1,\ldots, \MM_{t,t}\subseteq
\JJ_t$$
and
$$\RRR_{0,t},\ldots, \RRR_{t,t} \, .  $$
% In view of \eqref{def_rnn} and \eqref{def_rnm} this
%defines the collections $\RRR_{0,n},\ldots,\RRR_{n,n}$ and
%$\RRR_{n+1,n+1}$ as well.

We  now construct the corresponding  collections for $t=n+1$.

%\noindent {\em Stage 2.} Suppose that we have already constructed the collections
%$$\MM_{0,n}\subseteq \JJ_0, \   \MM_{1,n}\subseteq \JJ_1,\ldots, \MM_{n,n}\subseteq
%\JJ_n$$
%and
%$$\RRR_{0,n},\ldots, \RRR_{n,n} \, .  $$
%% In view of \eqref{def_rnn} and \eqref{def_rnm} this
%%defines the collections $\RRR_{0,n},\ldots,\RRR_{n,n}$ and
%%$\RRR_{n+1,n+1}$ as well.
%Now we will construct the collections
%$$
%\MM_{0,n+1}   \subseteq \JJ_0, \ \MM_{1,n+1}\subseteq \JJ_1,\ldots
%,\MM_{n+1,n+1}   \subseteq \JJ_{n+1}
%$$
%and
%$$\RRR_{0,n+1},\ldots,\RRR_{n+1,n+1}.$$

\vspace*{1ex}

\noindent{\em Stage 3.} Define
$$
\MM'_{n+1,n+1}:=\{J_{n+1}\in \JJ_{n+1}\;:\; J_{n+1}\subset J_n\mbox{ for some }J_n\in
\MM_{n,n}\}
$$
and let $\RRR_{n+1,n+1}$ be given by \eqref{def_rnn} with $m =n+1$.
Thus the collection  $\MM'_{n+1,n+1}$ consists of `good' intervals
from $\JJ_{n+1}$ that are contained within  some  interval from
$\MM_{n,n}$.  Our immediate task is to construct the corresponding
collections $ \MM'_{u,n+1}$ for each $ 0 \le u \le n$.  These will
be constructed  together with the `complementary'  collections
$\RRR_{u,n+1}$ in descending order with respect to $u$.

\vspace*{1ex}

\noindent{\em Stage 4.} With reference to Stage 3, suppose we have constructed the collections $\MM'_{u+1,n+1}$ and $\RRR_{u+1,n+1}$ for some $0\le u\le n$. We now construct  $\MM'_{u,n+1}$ and $\RRR_{u,n+1}$. Consider the collections $\MM_{u,n}$ and
$\RRR_{u,n}$. Observe that some of the intervals $J_u$ from $\MM_{u,n}$
may contain less than $R-[2R^{1-\epsilon/2}]$ sub-intervals from
$\MM'_{u+1,n+1}$ (or in other words, at least $[2R^{1-\epsilon/2}]$
intervals from $\RRR_{u+1,n+1}$). Such intervals $J_u$ fail the counting condition (C3)  for  $\MM_{u,n+1}$ and informally speaking are moved out of $\MM_{u,n}$ and into $\RRR_{u,n}$.  The resulting  sub-collections  are  $\MM'_{u,n+1}$ and $\RRR_{u,n+1}$ respectively. Formally,
$$
\MM'_{u,n+1}:=\{J_u\in \MM_{u,n}\;:\;
\#\{J_{u+1}\in \RRR_{u+1,n+1}\;:\; J_{u+1}\subset J_u\}<
[2R^{1-\epsilon/2}]\ \} \
$$
and $\RRR_{u,n+1}$ is given by \eqref{def_rnm} with $n=u$ and $ m=n+1$.

\vspace*{1ex}

\noindent{\em Stage 5.}
By construction the collections $\MM'_{u,n+1}$ satisfy conditions (C1) and (C3).  However,  for some $J_{u+1}\in \MM'_{u+1,n+1}$ it  may be the case that $J_{u+1}$ is not contained in any interval $J_u\in \MM'_{u,n+1}$  and thus the collections $\MM'_{u,n+1}$ are not necessarily nested. The point is that during  Stage 4 above the interval $J_u\in \JJ_u$ containing $J_{u+1}$ may be  `moved' into $\RRR_{u,n+1}$. In order to guarantee the nested condition (C2)  such intervals $J_{u+1}$ are removed from $\MM'_{u+1,n+1}$. The resulting sub-collection is the required auxiliary collection $\MM_{u+1,n+1}$. Note that $\MM_{u+1,n+1}$ is  constructed via $\MM'_{u+1,n+1}$ in ascending order with respect to $u$.  Formally,
$$
\MM_{0,n+1}:=\MM'_{0,n+1}
$$
and for $1\le u\le n+1$
$$
\MM_{u,n+1}:=\{J_u\in \MM'_{u,n+1}\;:\; J_u\subset J_{u-1} \mbox{
for some } J_{u-1}\in \MM_{u-1,n+1}\}   \, .
$$
With reference to Stage 2, this completes the  induction step and thereby   the construction of the auxiliary collections.

\vspace*{3ex}

For any integer $n \geq 0$, it remains to  construct the sought after collection   $\MM_n$ via the auxiliary collections $\MM_{n,m}$. Observe that since
$$
\MM_{n,n}\supset \MM_{n,n+1}\supset \MM_{n,n+2}\supset \ldots
$$
and the cardinality of each collection $\MM_{nm}$ with    $n\le m$ is finite, there exists some  integer $N(n)$ such that
$$
\MM_{n,m}   \ = \ \MM_{n,m'}  \qquad \forall \quad m,m'  \ge N(n)   \ .
$$

\noindent  Now simply define $$ \MM_n:=\MM_{n,N(n)} \ . $$

%\noindent There remains one slight issue with the above.  Namely that  $\MM_n$ could be empty.

\noindent Unfortunately, there remains one slight issue.  The collection   $\MM_n$  defined in this manner could be empty.

\vspace*{2ex}

The goal now  is to show  that  $\MM_{n,m}\neq
\emptyset$ for any $n\le m$. This clearly  implies  that $\MM_n\neq\emptyset$ and thereby completes the construction.

%\subsubsection{The collection $\MM_{m,n}$ is non-empty}

\begin{proposition}
 For all integers satisfying $ 0 \le n\le m$, the collection  $\MM_{n,m}$ is nonempty.
\end{proposition}

\noindent{\em Proof. \, } Suppose on the contrary that
$\MM_{n,m}=\emptyset$ for some integers satisfying $ 0 \le n\le m$.
In view of the construction of $\MM_{n,m}$,  every interval  from  $\MM_{n-1,m}$
contains at least $R-2R^{1-\epsilon/2}$ sub-intervals from $\MM_{n,m}$.
Therefore $\MM_{0,m}$ is empty and it follows that
$\RRR_{0,m}=\JJ_0$.

Now consider the set $\RRR_{n,m}$.  Note that
$$
\RRR_{n,m}\supseteq \RRR_{n,m-1}\supseteq\cdots\supseteq\RRR_{n,n}.
$$
and  that in view of \eqref{def_rnn}  elements of $\RRR_{n,n}$ are
intervals from $\II_n\backslash \JJ_n$.  Consider any interval
$J_n\in \RRR_{n,m}\backslash \RRR_{n,n}$. Then there exists an
integer $m_0$ with $n <  m_0 \le m $ such that $J_n\in \RRR_{n,m_0}$
but $J_n\not\in \RRR_{n,m_0-1}$. In view of \eqref{def_rnm}  any
interval  from $\RRR_{n,m_0}$ contains at least
$[2R^{1-\epsilon/2}]$ sub-intervals from $\RRR_{n+1,m_0}$ and
therefore from $\RRR_{n+1,m}$. The upshot is that for any  interval
$I_n\in \RRR_{n,m}$ we either have that $I_n\in \II_n\backslash
\JJ_n$ or that $I_n$ contains at least $[2R^{1-\epsilon/2}]$
intervals $I_{n+1}\in \RRR_{n+1,m}$.

%
%Note that the following condition about $\RRR_{mn}$ is true. Any
%segment $J_n\in \RRR_{mn}$ is either removed by the procedure from
%main theorem or contains at least $2R^{1-\epsilon/2}$ subsegments
%from $\RRR_{m+1,n}$.

Next we exploit Lemma \ref{sl_lem1}.
Choose an arbitrary interval $J_0$ from $ \RRR_{0,m} = \JJ_0$   and define
$\TTT_0:=\{J_0\}$. For $ 0\le n<m$, we define inductively the nested collections
$$
\TTT_{n+1}:=\{I_{n+1}\in \TTT(I_n)\;:\; I_n\in \TTT_n\} %\qquad 0\le n<m
$$
with $\TTT(I_n)$ given  by one of the following three scenarios.
\begin{itemize}
\item $I_n\in \RRR_{n,m}$ and $I_n$ contains at least $[2R^{1-\epsilon/2}]$ sub-intervals $I_{n+1}$ from $\RRR_{n+1,m}$. Let $\TTT(I_n)$ be any collection consisting of $[2R^{1-\epsilon/2}]$ such sub-intervals. Note that when $n=m-1$ we have $\TTT(I_n)\subset \RRR_{m,m}\subset \II_m\backslash \JJ_m$. Therefore $\TTT(I_{m-1})\cap \JJ_{m}=\emptyset$.

% Note that this case can only happen for $n<m$.
\item $I_n\in \RRR_{n,m}$ and $I_n$ contains strictly less than $[2R^{1-\epsilon/2}]$ sub-intervals $I_{n+1}$ from
$\RRR_{n+1,m}$. Then the interval $I_n \in \II_n\backslash \JJ_n$ and we  subdivide $I_n$ into $R$ closed intervals $I_{n+1}$ of equal length.  Let $\TTT(I_n)$ be any collection consisting of $[2R^{1-\epsilon/2}]$ such sub-intervals.  Note that $\TTT(I_n)\cap \JJ_{n+1}=~\emptyset$.

\item $I_n\not\in \RRR_{n,m}$. Then the interval $I_n$ does not
intersect any interval from $\JJ_n$ and we  subdivide $I_n$ into $R$ closed intervals $I_{n+1}$ of equal length.  Let $\TTT(I_n)$ be any collection consisting of $[2R^{1-\epsilon/2}]$ such sub-intervals. Note that $\TTT(I_n)\cap \JJ_{n+1}=~\emptyset$.
\end{itemize}

\noindent The upshot is that
$$
\# \TTT_n   \ =  \    \# \TTT_{n-1}  \times  [2R^{1-\epsilon/2}]   \qquad \forall  \quad  0< n\le m
$$
and that
$$
\TTT_{m}\cap \JJ_{m} \, =  \, \emptyset  \; .
$$
However, in view of Lemma~\ref{sl_lem1} the latter is impossible and therefore the starting premise that $\MM_{n,m}=\emptyset$ is false. This completes the proof of the proposition.
\newline\hspace*{\fill}$\boxtimes$

\newpage
\vspace*{8ex}

\noindent{\Large \bf Appendix:  The dual and  simultaneous forms of $\bad(i,j)$}

%\begin{center}
%\noindent{\Large \bf Appendix:  The dual and  `simultaneous' forms of $\bad(i,j)$}
%\end{center}

%\section{Appendix: The `dual $\equiv$  simultaneous' transference principle}
\vspace*{3ex}

\noindent Given a pair of real numbers $i$ and $j$ satisfying \eqref{neq1sv}, the following statement allows us to  deduce that the dual and simultaneous forms of $\bad(i,j)$ are equivalent.

\begin{theorem}\label{app_transference}
Let
$$
L_t(\vq) \, :=  \, \sum_s \theta_{ts} \, q_s   \qquad  (1\le s\le m, \, 1\le t\le n)
$$
be $n$ linear forms in $m$ variables and let
$$
M_s(\vu)\, :=  \, \sum_t \theta_{ts} \, u_t
$$
be the transposed set of $m$ linear forms in $n$ variables.  Suppose that there are integers
$\vq\neq \mathbf{0}$ such that
$$
||L_t(\vq)|| \, \le  \,  C_t  \, ,  \qquad   |q_s|\le X_s  \, ,
$$

\noindent for some constants $C_t$ and $X_s$ satisfying
$$
\max_s \{ \, D_s  :=  (l-1) \, X_s^{-1} \, d^{1/(l-1)}   \, \}    \, < \, 1   \
$$
where
$$ d  \, :=  \, \prod_{t} C_t  \, \prod_{s}
X_s  \,  \quad { and }  \quad   l\, :=  \, m+n \, .  $$
Then there are integers $\vu\neq
\mathbf{0}$ such that
\begin{equation}\label{app_th1_eq1}
||M_s(\vu)||  \, \le   \, D_s  \, ,  \qquad |u_t|\le U_t  \, ,
\end{equation}
where
\begin{equation*}%\label{app_th1_eq2}
U_t  \, :=  \, (l-1)  \, C_t^{-1}  \, d^{1/(l-1)}  \, .
\end{equation*}
\end{theorem}

%\begin{theorem}\label{app_transference}
%Let
%$$
%L_t(\vq) \, :=  \, \sum_s \theta_{ts} \, q_s   \qquad  (1\le s\le m, \, 1\le t\le n)
%$$
%be $n$ linear forms in $m$ variables and let
%$$
%M_s(\vu)\, :=  \, \sum_t \theta_{ts} \, u_t
%$$
%be the transposed set of $m$ linear forms in $n$ variables.  Suppose that there are integers
%$\vq\neq \mathbf{0}$ such that
%$$
%||L_t(\vq)|| \, \le  \,  C_t  \, ,\qquad |q_s|\le X_s  \, ,
%$$
%for some constants $C_t$ and $X_s$ satisfying  $0<C_t <1 \leq X_s$.
%Then there are integers $\vu\neq
%\mathbf{0}$ such that
%\begin{equation}\label{app_th1_eq1}
%||M_s(\vu)||  \, \le   \, D_s  \, ,  \qquad |u_t|\le U_t  \, ,
%\end{equation}
%where
%\begin{equation*}%\label{app_th1_eq2}
%D_s  \, :=   \, (l-1) \, X_s^{-1} \, d^{1/(l-1)}  \, ,\qquad
%U_t  \, :=  \, (l-1)  \, C_t^{-1}  \, d^{1/(l-1)}  \, ,
%\end{equation*}
%
%$$ d  \, :=  \, \prod_{t}C_t  \, \prod_{s}
%X_s  \,  \qquad { and }  \qquad   l\, :=  \, m+n \, .  $$
%\end{theorem}

This theorem is essentially a generalization of Theorem II  in \cite[ChapterV]{cassels}. In  short, compared to the latter,   the above theorem  allows the upper bounds for $||L_t(\vq)|| $   and $ |q_s| $ to vary with $ t$ and $s$ respectively.
The proof of Theorem~\ref{app_transference} makes use of  the following result  which  appears as Theorem~I in \cite[Chapter V]{cassels}.

\begin{proposition}\label{app_prop}
Let $f_k(\vz)$ ($1\le k\le l$) be $l$ linearly independent
homogeneous linear forms in the $l$ variables $\vz=(z_1,\ldots,z_l)$
and let $g_k(\vw)$ be $l$ linearly independent homogeneous linear
forms in the $l$ variables $\vw=(w_1,\ldots, w_l)$ of determinant
$d$. Suppose that all the products $z_iw_j$ ($1\le i,j\le l$) have
integer coefficients in
$$
\Phi(\vz,\vw)  \, :=  \, \sum_k f_k(\vz)  \, g_k(\vw)  \, .
$$
If the inequalities
$$
|f_k(\vz)|  \, \le  \,  \lambda \qquad (1\le k\le l)
$$
are soluble with integral $\vz\neq \mathbf{0}$ then the inequalities
$$
|g_k(\vw)|  \, \le \, (l-1) \, |\lambda \, d|^{1/(l-1)}  \, ,
$$
are soluble with integral $\vw\neq \mathbf{0}$.
\end{proposition}

\vspace*{2ex}

Armed with this proposition,  the proof of Theorem~\ref{app_transference} is relatively straightforward. Indeed, apart from obvious modifications the proof  is essentially as in \cite{cassels}.

\vspace*{2ex}

\noindent{\em Proof of Theorem \ref{app_transference}.  } We start by introducing  the  new variables
$$
\vp = (p_1,\ldots, p_n)  \qquad {\rm and}  \qquad  \vv =(v_1,\ldots, v_m)  \, .
$$
Now let

$$
f_k(\vq,\vp) \, :=  \, \left\{\begin{array}{lcl} C_k^{-1} \, (L_k(\vq)+p_k)  \  &\mbox{
if }&1\le k\le n\\[2ex]
X_{k-n}^{-1} \, q_{k-n}  \ &\mbox{ if }&n<k\le l
\end{array}\right.
$$
and

$$
g_k(\vu,\vv)  \, :=   \, \left\{\begin{array}{lcl} C_k u_k&\mbox{
if }&1\le k\le n\\[2ex]
X_{k-n}(-M_{k-n}(\vu)+v_{k-n})&\mbox{ if }&n<k\le l   \, .
\end{array}\right.
$$

\noindent Then the $f_k$ are linearly independent forms in the $l:=m+n$ variables
$\vz=(\vq,\vp)$ and the $g_k$ are linearly independent forms in the
$l$ variables $\vw=(\vu,\vv)$ with determinant
$$
d   \, :=  \, \prod_{t=1}^nC_t  \, \prod_{s=1}^m X_s \, .
$$
Furthermore,
$$
\sum_{k\le l} f_kg_k \, =   \, \sum_{t\le n}u_tp_t+\sum_{s\le m}v_sq_s   \,
$$
since the terms in $u_tq_s$ all cancel out. By hypothesis there are
integers $\vq \neq\mathbf{0} $ and $ \vp$ such that
$$
|f_k(\vq,\vp)|\le 1,
$$
so we may apply Proposition~\ref{app_prop} with $\lambda=1$. It follows that there are integers $(\vu,\vv)\neq(\mathbf{0},\mathbf{0})$ such that
\begin{equation*}  \label{yippy}
\left. \begin{array}{l} C_t|u_t|\\[1ex]
X_s|-M_s(\vu)+v_s|
\end{array}\right\}  \ \le \  (l-1) \, d^{1/(l-1)}   \;
\end{equation*}
and so the inequalities given by \eqref{app_th1_eq1} hold.  It remains to show that $\vu\neq\mathbf{0}$. By hypothesis $D_s < 1 $ for all $s$  and so if  $ \vu=\mathbf{0}$ we must have that $v_s = 0 $ for all $s$.      However   $(\vu,\vv)=(\mathbf{0},\mathbf{0})$ is excluded.
\newline\hspace*{\fill}$\boxtimes$

\vspace*{4ex}

\noindent Given Theorem~\ref{app_transference}, it  is relatively straightforward to show that the  dual and simultaneous forms of $\bad(i,j)$ are equivalent.

\vspace*{2ex}

 Suppose the point $(x,y)\in \RR^2$ does not belong to the
simultaneous $\bad(i,j)$ set. It follows from the definition of the latter that for any
constant $c>0$ there exists an integer $q_0 \geq 1 $ such that
$$
 \begin{array}{l}
||q_0x||  \, \le \,  c \,  q_0^{-i}\\[1ex]
||q_0y||  \, \le  \,  c \, q_0^{-j}  \, .
\end{array}
$$
Without loss of generality assume that $c < 1/2$. With reference to Theorem~\ref{app_transference},  let  $m=1$, $n=2$, $L_1(\vq)=qx$, $L_2(\vq)=qy$,
$C_1=cq_0^{-i}$,  $ C_2=cq_0^{-j}  $  and $ X_1=q_0$. Hence there exists an  integer pair
$(u_1,u_2) \neq (0,0)$ such that
$$
\begin{array}{l}
|| xu_1+yu_2 ||  \, \le  \,  2cq_0^{-1}\\[1ex]
|u_1|  \, \le   \, 2q_0^i\\[1ex]
|u_2| \, \le  \,  2q_0^j   \; .
\end{array}
$$
This in turn implies that
\begin{equation} \label{slv10}
\max\{|u_1|^{1/i}, |u_2|^{1/j}\}  \ ||xu_1+yu_2||  \; \le  \;
2^{1/i+1/j+1}  \, c  \; .
\end{equation}
In other words, for any  arbitrary small constant $c > 0$ there exists $(u_1,u_2)  \in\ZZ^2 \backslash \{(0,0)\}$ for which  \eqref{slv10} is satisfied.  It follows that the point  $(x,y)$ does not belong to the dual
$\bad(i,j)$ set.  The upshot  is that the
dual  $\bad(i,j)$ set is
a subset of the simultaneous $\bad(i,j)$ set.

\vspace*{1ex}

Suppose the point $(x,y)\in \RR^2$ does not belong to the
dual $\bad(i,j)$ set. It follows from the definition of the latter that for any
constant $c>0$ there exists $(a,b)  \in\ZZ^2 \backslash \{(0,0)\}$ such that
$$
\max\{|a|^{1/i}, \,  |b|^{1/j}\}  \ ||ax+by|| \, \le \,  c.
$$
Without loss of generality assume that $c < 1/4$ and let $ q_0 :=  \max\{|a|^{1/i}, \, |b|^{1/j}\} $. With reference to Theorem~\ref{app_transference},  let  $m=2$,  $n=1$, $L_1(\vq)=q_1x+q_2y$,  $ C_1=cq_0^{-1}$,  $ X_1=q_0^i$  and  $X_2=q_0^j$.    Hence there exists an  integer
$u\neq 0$ such that
$$
\begin{array}{l}
||u x||  \, \le \, 2 \, c^{1/2}  \, q_0^{-i}\\[1ex]
||u y||  \, \le \, 2 \, c^{1/2}  \, q_0^{-j}\\[1ex]
|u|  \, \le  \,  2 \, c^{-1/2}  \, q_0  \ .
\end{array}
$$
This in turn implies that there exists and integer $q = | u| \ge 1$ such that
\begin{equation} \label{slv101}
\max\{||qx||^{1/i}, ||qy||^{1/j}\} \; \le \;
\max\left\{2^{\frac{1+i}{i}}c^{\frac{j}{2i}}, \,
2^{\frac{1+j}{j}}c^{\frac{i}{2j}}\right\} \,  q^{-1}   \ .
\end{equation}
In other words, for any  arbitrary small constant $c > 0$ there exists $q \in \NN $ for which  \eqref{slv101} is satisfied.  It follows that the point  $(x,y)$ does not belong to the simultaneous
$\bad(i,j)$ set.  The upshot  is that the
 simultaneous $\bad(i,j)$ set is
a subset of the dual $\bad(i,j)$ set.

\vspace*{6ex}

\noindent{\bf Acknowledgements.} SV would like to thank Graham Everest for being such a pillar of  support throughout his mathematical life -- especially during the teenage years!  As ever, an enormous thankyou to Bridget, Ayesha and Iona for just about everything so far this millennium.

%\newpage

\def\cprime{$'$}

\vspace{10mm}

\noindent Dzmitry A. Badziahin: Department of Mathematics,
University of York,

\vspace{0mm}

\noindent\phantom{Dzmitry A. Badziahin: }Heslington, York, YO10 5DD,
England.

%\vspace{0mm}

\noindent\phantom{Dzmitry A. Badziahin: }e-mail: db528@york.ac.uk

\vspace{5mm}

\noindent Andrew D. Pollington: National Science Foundation

\vspace{0mm}

\noindent\phantom{Andrew D. Pollington: }Arlington VA 22230 USA

%\vspace{0mm}

\noindent\phantom{Andrew D. Pollington: }e-mail: adpollin@nsf.gov

\vspace{5mm}

\noindent Sanju L. Velani: Department of Mathematics, University of
York,

\vspace{0mm}

\noindent\phantom{Sanju L. Velani: }Heslington, York, YO10 5DD,
England.

%\vspace{0mm}

\noindent\phantom{Sanju L. Velani: }e-mail: slv3@york.ac.uk


\begin{thebibliography}{10}


\bibitem{cassels}
J.W.S.~Cassels: {\em An introduction to the geometry of numbers.}
  Classics in Mathematics, Springer-Verlag, Berlin, (1997).
\newblock Corrected reprint of the 1971 edition.


\bibitem{dav}
H.~Davenport: A note on Diophantine approximation II. {\em
Mathematika} \textbf{11} (1964) 50--58.


\bibitem{falc}
K.~Falconer: {\em Fractal Geometry: Mathematical Foundations and
Applications.} John Wiley \& Sons, (1990).


%\bibitem{ekl}
%M.~Einsiedler, A.~Katok and E.~Lindenstrauss :  Invariant measures
%and the set of exceptions to Littlewood's conjecture.  {\em Annals
%of Math.} \textbf{164}  (2006),  513--560.

\bibitem{KW}
D.~Kleinbock and B.~Weiss:  Modified Schmidt games and Diophantine
approximation with weights.  {\em Adv. Math.} to appear. Pre-print: arXiv:0805.2934 (2008), 1--22.


\bibitem{KTV}
S.~Kristensen, R.~Thorn and S.L.~Velani: Diophantine approximation
and badly approximable sets. {\em Adv. Math.} \textbf{203} (2006),
132--169.


\bibitem{pvl}
A.D.~Pollington and  S.L.~Velani: On a problem in simultaneously
Diophantine approximation: Littlewood's conjecture. {\em Acta
Math.} \textbf{66} (2000),  29--40.


\bibitem{PV}
A.D.~Pollington and  S.L.~Velani: On simultaneously badly
approximable pairs. {\em Jou. Lond. Math. Soc.} \textbf{66} (2002),
 29--40.



\bibitem{schconj}
W.M.~Schmidt: Open problems in Diophantine approximation. {\em
Approximations diophantiennes et nombres transcendants (Luminy
1982)}, Progress in Mathematics, Birkh\"auser, (1983).

\bibitem{vent}
A.~Venkatesh: The work of Einsiedler, Katok and Lindenstrauss on
the Littlewood conjecture.  {\em Bull. Amer. Math. Soc.} \textbf{45}
(2008), 117--134.



\end{thebibliography}
\end{document}